 \documentclass[10pt,twoside,english]{article}
\usepackage{amssymb,amsmath,babel,geometry,latexsym,graphics,tabularx,shapepar,enumerate}
\ifx\optionkeymacros\undefined\else \fi

\catcode`\Œ=\active\defŒ{{\aa}}       
\catcode`\º=\active\defº{\int}        
\catcode`\=\active\def{\c c}        
\catcode`\¶=\active\def¶{\partial}    
\catcode`\Ä=\active\defÄ{\oint}       
\catcode`\Æ=\active\defÆ{\triangle}   
\catcode`\Â=\active\defÂ{\neg}        
\catcode`\µ=\active\defµ{\mu}         
\catcode`\¿=\active\def¿{{\o}}        
\catcode`\¹=\active\def¹{\pi}         
\catcode`\Ï=\active\defÏ{{\oe}}       
\catcode`\§=\active\def§{{\ss}}       
\catcode`\ =\active\def {\dagger}     
\catcode`\Ã=\active\defÃ{\sqrt}       
\catcode`\·=\active\def·{\Sigma}      
\catcode`\Å=\active\defÅ{\approx}     
\catcode`\½=\active\def½{\Omega}      
\catcode`\£=\active\def£{{\it\$}}     
\catcode`\°=\active\def°{\infty}      
\catcode`\¤=\active\def¤{{\S}}        
\catcode`\¦=\active\def¦{{\P}}        
\catcode`\¥=\active\def¥{\bullet}     
\catcode`\»=\active\def»{\leavevmode\raise.585ex\hbox{\b a}}      
\catcode`\¼=\active\def¼{\leavevmode\raise.6ex\hbox{\b o}}        
\catcode`\­=\active\def­{\not=}       
\catcode`\²=\active\def²{\leq}        
\catcode`\³=\active\def³{\geq}        
\catcode`\Ö=\active\defÖ{\div}        
\catcode`\É=\active\defÉ{{\dots}}     
\catcode`\¾=\active\def¾{{\ae}}       
\catcode`\Ç=\active\defÇ{\ll}         
\catcode`\Ò=\active\defÒ{``}          
\catcode`\Á=\active\defÁ{!`}          
\catcode`\¢=\active\def¢{\rlap/c}     
\catcode`\Ô=\active\defÔ{`}           
\catcode`\Õ=\active\defÕ{'}           

\catcode`\=\active\def{{\AA}}       
\catcode`\'=\active\def'{\c C}        
\catcode`\¯=\active\def¯{{\O}}        
\catcode`\¸=\active\def¸{\Pi}         
\catcode`\Î=\active\defÎ{{\OE}}       
\catcode`\®=\active\def®{{\AE}}       
\catcode`\×=\active\def×{\diamond}    
\catcode`\¡=\active\def¡{\accent'27}  
\catcode`\Ó=\active\defÓ{''}          
\catcode`\±=\active\def±{\pm}         
\catcode`\È=\active\defÈ{\gg}         
\catcode`\À=\active\defÀ{?`}          
\catcode`\Ð=\active\defÐ{--}          
\catcode`\Ñ=\active\defÑ{---}         

\catcode`\Š=\active\defŠ{\"a}        
\catcode`\'=\active\def'{\"e}        
\catcode`\•=\active\def•{\"{\i}}     
\catcode`\š=\active\defš{\"o}        
\catcode`\Ÿ=\active\defŸ{\"u}        
\catcode`\Ø=\active\defØ{\"y}        
\catcode`\€=\active\def€{\"A}        
\catcode`\…=\active\def…{\"O}        
\catcode`\†=\active\def†{\"U}        
\catcode`\‡=\active\def‡{\'a}        
\catcode`\Ž=\active\defŽ{\'e}        
\catcode`\'=\active\def'{\'{\i}}     
\catcode`\—=\active\def—{\'o}        
\catcode`\œ=\active\defœ{\'u}        
\catcode`\ƒ=\active\defƒ{\'E}        
\catcode`\ˆ=\active\defˆ{\`a}        
\catcode`\=\active\def{\`e}        
\catcode`\"=\active\def"{\`{\i}}     
\catcode`\˜=\active\def˜{\`o}        
\catcode`\=\active\def{\`u}        
\catcode`\Ë=\active\defË{\`A}        
\catcode`\‹=\active\def‹{\~a}        
\catcode`\–=\active\def–{\~n}        
\catcode`\›=\active\def›{\~o}        
\catcode`\Ì=\active\defÌ{\~A}        
\catcode`\"=\active\def"{\~N}        
\catcode`\Í=\active\defÍ{\~O}        
\catcode`\‰=\active\def‰{\^a}        
\catcode`\=\active\def{\^e}        
\catcode`\"=\active\def"{\^{\i}}     
\catcode`\™=\active\def™{\^o}        
\catcode`\ž=\active\defž{\^u}        

\let\optionkeymacros\null

\usepackage{hyperref}
\usepackage[all,2cell]{xy} \UseAllTwocells \SilentMatrices
\def\carlitz{{{\text{Car}}}}

\newtheorem{Theorem}{Theorem}
\newtheorem{Lemme}[Theorem]{Lemma}
\newtheorem{Proposition}[Theorem]{Proposition}
\newtheorem{Corollaire}[Theorem]{Corollary}
\newtheorem{Definition}[Theorem]{Definition}

\newcommand{\ol}[1]{\chi_t(#1)}

\newcommand{\ZZ}{\mathbb{Z}}
\newcommand{\FF}{\mathbb{F}}
\newcommand{\CC}{\mathbb{C}}
\newcommand{\QQ}{\mathbb{Q}}
\newcommand{\RR}{\mathbb{R}}

\newcommand{\GG}{\mathbb{G}}
\newcommand{\KK}{\mathbb{K}}
\newcommand{\LL}{\mathbb{L}}
\newcommand{\MM}{\mathbb{M}}

\newcommand{\TT}{\mathbb{T}}

\newcommand{\bfs}{\boldsymbol{s}}

\newcommand{\bsb}{\boldsymbol}

\newcommand{\sqm}[4]{
\left(\begin{array}{ll}#1 & #2 \\ #3 & #4\end{array}\right)}
 
\newcommand\CVD{{\hfill\hfil{\lower 2 pt\hbox{\vrule\vbox to 7pt 
{\hrule width 6pt\vfill\hrule}\vrule}}}\vskip 0.5cm}

\def\SL{{\bf SL}}
\def\GL{{\bf GL}}

\let\ge=\geqslant                                 

\title{$\tau$-recurrent sequences and modular forms}
\author{Federico Pellarin\footnote{Current address: LaMUSE, 23, rue du Dr. Paul Michelon, 42023 Saint-Etienne Cedex.} \footnote{Supported by the contract ANR ``HAMOT", BLAN-0115-01.}}

\begin{document}

\maketitle

\begin{small}
\noindent\textbf{Abstract.} In this paper we deal with {\em Drinfeld
modular forms}, defined and taking values in complete fields of positive characteristic. Our aim is to study 
a sequence $(g_k^\star(z,t))_{k\geq 0}$ of 
families of Drinfeld modular forms that produces, for certain values of the parameter $t$, 
several kinds of Eisenstein series considered by Gekeler. We obtain
formulas involving these functions depending on the parameter $t$. To obtain our results, we 
introduce and discuss {\em $\tau$-linear recurrent sequences} and {\em deformations of vectorial
modular forms}, and we give on the way some applications to some special values of $L$-functions and to the study of {\em extremal quasi-modular forms}.
\end{small}


%

\medskip

\setcounter{tocdepth}{1}
\tableofcontents

\section{Introduction, results}

The present paper deals with the following loosely question: {\em why do there exist two kinds of ``Eisenstein series" in the theory of Drinfeld 
modular forms for $\mathbf{GL}_2(\FF_q[\theta])$, and not just one as in the theory for $\mathbf{SL}_2(\ZZ)$?} Here, we give a tentative of answer: {\em the theories are not 
too different, the two kinds of the drinfeldian framework indeed come from a unique family of functions.}
We will study this family and give some applications of our investigations.

\medskip

let $q=p^e$ be a power of a
prime number $p$ with $e>0$ an integer, let $\FF_q$ be the finite field with $q$ elements.
We consider the polynomial ring $A=\FF_q[\theta]$ and its fraction field $K=\FF_q(\theta)$, with $\theta$
an indeterminate over $\FF_q$. On $K$, we will consider the absolute value
$|\cdot|$ defined by $|a|=q^{\deg_\theta a}$, $a$ being in $K$, so that
$|\theta| = q$.  Let $K_\infty :=\FF_q((1/\theta))$ be the
completion of $K$ for this absolute value, let $K_\infty^{\text{alg.}}$ be an
algebraic closure of $K_\infty$, let $\CC_\infty$ be the completion of
$K_\infty^{\text{alg.}}$ for the unique extension of $|\cdot|$ to $K_\infty^{\text{alg.}}$, and let $K^{\text{alg.}}$ be the algebraic closure of $K$
in $\CC_\infty$. The presentation of our results requires that we first introduce some of the tools that will be used all along the paper.

\medskip\noindent\emph{1. Drinfeld modular and quasi-modular forms}.
Following Gekeler in \cite{Ge}, we denote by $\Omega$ the set $\CC_\infty\setminus K_\infty$, which has a structure of 
rigid analytic space.
The group $$\Gamma=\mathbf{GL}_2(A)$$ acts discontinuously on $\Omega$ by homographies; for $\gamma=\begin{pmatrix}a&b\\ c&d\end{pmatrix}\in\Gamma$
and $z\in \Omega$,
we denote by $\gamma(z)=(az+b)/(cz+d)$ the action of $\gamma$ on $z$. 
Gekeler considered three algebraically independent functions $$E,g,h:\Omega\rightarrow \CC_\infty$$
such that, for all $\gamma=\begin{pmatrix}a&b\\ c&d\end{pmatrix}\in\Gamma$ and $z\in\Omega$:
\begin{eqnarray}
g(\gamma(z))&=&(cz+d)^{q-1}g(z),\nonumber\\
h(\gamma(z))&=&(cz+d)^{q+1}\det(\gamma)^{-1}h(z),\nonumber\\
E(\gamma(z))&=&(cz+d)^2\det(\gamma)^{-1}\left(E(z)-\frac{c}{\widetilde 
{\pi}(cz+d)}\right).\label{formE}
\end{eqnarray}
These functions are holomorphic in the sense of \cite[Definition 2.2.1]{FP}.
Here, $\widetilde{\pi}$ is a fundamental period of the {\em Carlitz exponential function} $e_{\carlitz}$ defined, for all $\zeta\in \CC_\infty$, by the sum of the converging series:
\begin{equation}\label{exponential}
e_\carlitz(\zeta)=\sum_{n\ge 0}\frac{\zeta^{q^n}}{d_n},
\end{equation}
where $d_0:=1$ and $d_i:=[i][i-1]^q\cdots[1]^{q^{i-1}}$, with $[i]=\theta^{q^i}-\theta$ if $i>0$.

It is possible to show that $\widetilde{\pi}$ is equal, up to a choice of a $(q-1)$-th root of $-\theta$,
to the (value of the) convergent product:
$$\widetilde{\pi}:=\theta(-\theta)^{\frac{1}{q-1}}\prod_{i=1}^\infty(1-\theta^{1-q^i})^{-1}\in K_\infty((-\theta)^{\frac{1}{q-1}})\setminus K_\infty.$$

According to Gekeler in \cite{Ge} (but we will prefer to borrow notations from Gerritzen-van der Put in \cite{GePu}), the ``local parameter at infinity" of the quotient space 
$\Gamma\backslash\Omega$ can be defined as a map $\Omega\rightarrow \CC_\infty$ by:
$$u(z)=\frac{1}{e_\carlitz(\widetilde{\pi}z)}.$$ In \cite{Ge}, it is proved that $E,g,h$ have, locally at $u=u(z)=0$, convergent $u$-expansions in $A[[u]]$.
The functional equations above and the ``nice local behavior at infinity" indicate that $g,h$ are {\em Drinfeld modular forms},
of {\em weights} $q-1$, $q+1$ and {\em types} $0,1$ respectively. 

After (\ref{formE}) it is apparent that $E$ is not a Drinfeld modular form. In \cite{Ge},
Gekeler calls it a ``false Eisenstein series" of weight $2$ and type $1$. Nevertheless, it is the prototype of {\em Drinfeld quasi-modular form}, of weight 
$2$, type $1$ and depth $1$ (see \cite{BP}).

\medskip

In the classical theory of modular forms for $\mathbf{SL}_2(\ZZ)$, Eisenstein series
\begin{equation}\label{defeisensteinseries}G_{2\nu}(z)=\sideset{}{'}\sum_{c,d\in \ZZ}(cz+d)^{-2\nu}=2\zeta(2\nu)+2\frac{(2\pi \nu)^{2\nu}}{(2\nu-1)!}\sum_{n=1}^\infty\sigma_{2\nu-1}(n)q^n,\quad \nu\geq 2\end{equation}
appear in the coefficients of the Laurent expansion 
at zero of Weierstrass $\wp$-functions
\begin{equation}\label{weierstrass}\wp_{z\ZZ+\ZZ}(\zeta)=\zeta^{-2}+\sum_{\nu=1}^\infty(2\nu+1)G_{2\nu+2}(z)\zeta^{2\nu},
\end{equation} whose differential equation yields explicit quadratic recursive relations
\begin{equation}\label{quadraticrec}
(2\nu-2)(2\nu+5)F_{2\nu+4}=12(F_4F_{2\nu}+F_6F_{2\nu-2}+\cdots+F_{2\nu}F_4),
\end{equation}
where $F_{2\nu}=(-1)^\nu(2\nu-1)2^{-1-2\nu}\pi^{-2\nu}G_{2\nu}$.
A proof can be found in Weil's book \cite[Chapter 5, (1)]{We} but the formula (\ref{quadraticrec}), in its precise form, is copied from \cite{SD}.

Goss was the first, in  \cite{Go1}, to show the existence of $u$-expansions for a variant of Eisenstein series 
associated with the group $\Gamma$:
\begin{equation}\label{defeisensteinseries2}G_{k(q-1)}=\sideset{}{'}\sum_{c,d\in A}(cz+d)^{-k(q-1)},\quad(k\geq1)
\end{equation} (the dash $'$ means that we avoid the couple $(c,d)=(0,0)$ in the sum, and we allow an abuse of notation here).
These series no longer occur in the series expansion of
the ``analogues of $\wp$": the {\em exponential functions} of rank $2$ {\em Drinfeld modules}. 

For $z\in\Omega$, we denote by $\Lambda_z$ the $A$-module $A+zA$, free of rank $2$.
The evaluation at $\zeta\in \CC_\infty$ of the exponential function $e_{\Lambda_z}$ associated to the lattice $\Lambda_z$ is given by the series 
\begin{equation}\label{ellipticexponential}
e_{\Lambda_z}(\zeta)=\sum_{i=0}^\infty\alpha_i(z)\zeta^{q^i},
\end{equation}
for  functions $\alpha_i:\Omega\rightarrow \CC_\infty$ with $\alpha_0=1$. This series expansion can be understood in many ways as
the analogue of the Laurent series expansion at $0$ for the Weierstrass $\wp$-function in the Drinfeldian theory. For all $i$, $\alpha_i$ is a Drinfeld modular forms of weight $q^i-1$ and 
type $0$, but in general, it is not proportional to $G_{q^i-1}$.

Following Gekeler in \cite[(2.8)]{Ge}, the generating series of the $G_{k(q-1)}$'s is indeed related to the inverse of the series $e_{\Lambda_z}$:
$$\frac{\zeta}{e_{\lambda_z}(\zeta)}=1-\sum_{k\geq 1}G_{k(q-1)}(z)\zeta^{k(q-1)}.$$
In contrast with the use one can make of (\ref{weierstrass}), no general recursive formula is known to be deducible from the above generating series, involving all the series $G_{k(q-1)}$, like (\ref{quadraticrec}).

At least, the latter generating series can be used to show that a zero density proportion of of the $G_{k(q-1)}$'s
occur in the series expansion of the composition inverse of $e_{\Lambda_z}$, the {\em logarithm} series $\log_{\Lambda_z}$ associated to $\lambda_z$:
\begin{equation}\label{ellipticexponential}
\log_z(\mu)=\sum_{k\geq 1}G_{q^k-1}(z)\mu^{q^k}.\end{equation}
The functions $e_{\Lambda_z}$ and $\log_{\Lambda_z}$ do not seem to satisfy interesting differential relations. Instead of this, they satisfy 
functional equations. Let 
$$\phi_z(\theta)=\theta\tau^0+\widetilde{\pi}^{q-1}g(z)\tau^1+\widetilde{\pi}^{q^2-1}\Delta(z)\tau^2,$$
in $\CC_\infty[\tau]=\mathbf{End}_{\FF_q-\text{lin.}}(\GG_a(\CC_\infty))$ be the elliptic Drinfeld module associated with $\Lambda_z$. 
We then have, for all $a\in A$,
\begin{equation}\label{functionalequations}\phi_z(a)e_z(\zeta)=e_z(a\zeta),\quad \zeta\in \CC_\infty,\quad \log_z(\phi_z(\theta)\mu)=\theta\log_z(\mu),\quad |\mu| \text{ small enough}.
\end{equation}
These functional equations yield recursive formulas involving the coefficients $\alpha_i,G_{q^k-1}$ of the respective series expansions
at $0$ of $e_{\Lambda_z}$ and $\log_{\lambda_z}$.

As Gekeler points out, for example in \cite{Gek2}, we must then distinguish between
{\em ortho-}Eisenstein (\footnote{We drop the suffix ``ortho" henceforth.}) series and {\em para-}Eisenstein series, belonging to sequences
that we recall now.

The first sequence, $(g_k)_{k\geq 0}$, introduced by Gekeler in \cite{Ge},
is determined by setting $g_0=1, g_1=g$
and then, inductively along the first formula of (\ref{functionalequations}), defining:
$$g_k=g_{k-1}g^{q^{k-1}}-[k-1]g_{k-2}
\Delta^{q^{k-2}},\quad (k\ge 2),$$
where $\Delta=-h^{q-1}$. 
For fixed $k$, this is the normalisation (\footnote{A formal Laurent series in powers of $u$ is said to be {\em normalised} if the monomial
of lowest order in $u$ appearing in it is monic. By abuse of language, we will say that a non-zero modular form is normalised if its $u$-expansion 
is normalised.}) of the Eisenstein series of weight $q^k-1$:
\begin{equation}\label{eq:definition_eisenstein}
g_k=(-1)^{k+1}\widetilde{\pi}^{1-q^k}[k][k-1]\cdots[1]G_{q^k-1}.
\end{equation}

Together with $(g_k)_{k\geq 0}$, we have the second sequence $(m_k)_{k\geq 0}$ of Drinfeld modular forms, called
{\em para-Eisenstein series}, discussed by Gekeler in \cite{Gek2}, possessing the 
same sequence of weights and types. This can be defined inductively along the second formula of (\ref{functionalequations}):
$$m_0=1,\quad m_1=g,\quad m_k=gm_{k-1}^q+[k-1]^q\Delta m_{k-2}^{q^2},\quad (k\geq 2),$$ as it is easily shown that $\alpha_k=d_km_k$
for all $k$. Several authors have studied the sequences of functions $(g_k)_{k\geq0},(m_k)_{k\geq 0}$ noticing several 
similarities and differences (congruences, location of zeros, behaviour under the action of Hecke's operators etc.).
We introduce now the sequence of functions studied here.

\medskip

\noindent\emph{The sequence $g_k^\star$.} For $w$ integer and $m\in\ZZ/(q-1)\ZZ$, we denote by $M_{w,m}$ the $\CC_\infty$-vector space 
of Drinfeld modular forms of weight $w$ and type $m$ and we denote by $M=\oplus_{w,m}M_{w,m}$ the $\CC_\infty$-algebra generated by Drinfeld 
modular forms. Gekeler \cite[Theorem (5.13)]{Ge} proved that $M=\CC_\infty[g,h]$. 

Let $$\tau:M[z]\otimes \FF_q((t))\rightarrow M[z]\otimes \FF_q((t))$$ be the unique $\FF_q((t))$-linear
map extending the Frobenius map $f\mapsto f^q$ on $M[z]$; for example, $\tau t=t,\tau g=g^q,\tau z=z^q$. It is easy to show that $\tau$ induces an injective linear map from $M_{w,m}\otimes \FF_q((t))$ to $M_{qw,m}\otimes \FF_q((t))$.

We introduce here a third new sequence $(g^\star_k)_{k\geq0}$ recursively as follows: we take once again the initial data $g^\star_0=1, g^\star_1=g$ and
then we set:
\begin{equation}\label{firstrecurrence}
g^\star_k=g(\tau g^\star_{k-1})+(t-\theta^q)\Delta(\tau^2g^\star_{k-2}),\quad (k\geq 2).
\end{equation}
For example,
\begin{eqnarray*}
g^\star_2&=& g^{1+q}-\Delta (t-\theta^q),\\
g^\star_3&=& g^{1+q+q^2}+\Delta g^{q^2}(t-\theta^q)+\Delta^q g(t-\theta^{q^2}), 
\end{eqnarray*}
and in general, $g^\star_k\in M_{q^k-1,0}\otimes\FF_q[t]$.
The recursion process involved in this definition is typical of what we study in this paper; we will refer to such a kind of sequence 
as to a {\em $\tau$-linear recurrent sequence} (of order $2$, because $g^\star_k=X(\tau g^\star_{k-1})+Y(\tau^2g^\star_{k-2})$
with $X=g$ and $Y=(t-\theta^q)\Delta$).

It is easy to show (by induction, but read Section \ref{recurrencesequences}) that the sequence $(g^\star_k)_{k\geq0}$ can be constructed alternatively, by 
choosing the same initial data and setting:
\begin{equation}\label{secondrecurrence}
g^\star_k=g^{q^{k-1}}g^\star_{k-1}+\Delta^{q^{k-2}}(t-\theta^{q^{k-1}})g^\star_{k-2},\quad (k\geq 2).\end{equation} 
This time, we speak about a {\em $\tau$-linearised recurrent sequence} (of order $2$, because we can write
$g^\star_k=(\tau^kX)g^\star_{k-1}+(\tau^kY)g^\star_{k-2}$ with $X=g^{1/q}$ and $Y=\Delta^{1/q^2}(t-\theta^{1/q})$). 
The rudiments of the theory of such sequences are established in
Section \ref{recurrencesequences} below.

These two expressions of the same sequence have particular benefit for us.
Indeed, from (\ref{firstrecurrence}), it follows that 
$$g^\star_k(z,\theta^{q^k})=m_k(z),\quad (k\geq 0),$$ while from
(\ref{secondrecurrence}) we obtain that
$$g^\star_k(z,\theta)=g_k(z),\quad (k\geq 0).$$

Summing up everything together, we understand that the sequence $(g^\star_k)_{k\geq 0}$ somewhat encompasses the
two different kinds of Eisenstein series of weights $(q^k-1)_{k\geq 0}$ considered by Gekeler and is an object that needs to be studied on its own. 
This we do in the present paper, but, before arriving at our principal results, we need to introduce further tools.

\medskip\noindent\emph{The functions $s_\carlitz,\bsb{d}_1,\bsb{d}_2$.} 
We need to recall the definitions of some now classical functions that naturally arise in the theory 
of Anderson's $t$-motives, also considered in \cite{archiv} (we will use notations introduced in the latter paper).
The first function is defined by the series:
\begin{equation}\label{scarlitz}
s_{\carlitz}(t):=\sum_{i=0}^\infty e_\carlitz\left(\frac{\widetilde{\pi}}{\theta^{i+1}}\right)t^i=\sum_{n=0}^\infty\frac{\widetilde{\pi}^{q^n}}{d_n(\theta^{q^n}-t)},
\end{equation} both converging for $|t|<q$. This is the {\em canonical rigid analytic trivialisation} of the so-called {\em Carlitz's motive}.
We refer to \cite{Bourbaki} for a description of the main properties of
it, or to the papers \cite{An, ABP,Pa}, where it was originally introduced and appears with different notations.

The series (introduced and studied in  \cite{archiv})
\begin{eqnarray*}
\bsb{s}_1(z,t)&=&\sum_{i=0}^\infty\frac{\alpha_i(z)z^{q^i}}{\theta^{q^i}-t},\\
\bsb{s}_2(z,t)&=&\sum_{i=0}^\infty\frac{\alpha_i(z)}{\theta^{q^{i}}-t},
\end{eqnarray*} converge on $\Omega\times B_{q}$ and define two
functions $\Omega\rightarrow \CC_\infty[[t]]$ with the series in the image converging on $B_{q}$
 ($B_r$ denotes the open disk of center $0$ and radius $r>0$). We point out that for a fixed choice of
$z\in\Omega$,
the matrix function $(\bsb{s}_1(z,t),\bsb{s}_2(z,t))$ is the {\em canonical rigid analytic trivialisation} of the {\em $t$-motive associated}
to the lattice $\Lambda_z$.
We recall that we have set, for $i=1,2$ (in \cite{archiv}):
$$\bsb{d}_i(z,t):=\widetilde{\pi}s_\carlitz(t)^{-1}\bfs_i(z,t)$$
and we point out that, in the notations of \cite{archiv}, $\bsb{d}_2=\bsb{d}$.
The advantage of using these functions, comparing with the $\bsb{s}_i$'s, is that evaluation at $t=\theta$ makes sense, and we can check:
\begin{equation}\label{limitsd}\bsb{d}_1(z,\theta)=z,\quad \bsb{d}_2(z,\theta)=1.\end{equation}
Moreover, as it was pointed out in \cite{archiv}, $\bsb{d}_2$ has a $u$-expansion defined over $\FF_q[t,\theta]$ (see later, Proposition \ref{prarchiv}).
We then have the following result.

\begin{Theorem}\label{charroots} For all $k\geq 0$, we have, for $(z,t)\in\Omega\times \CC_\infty$ with $|t|<q^{k+1}$ and writing $g_k^{\star}=g_k^{\star}(z,t)$,
$s_\carlitz=s_\carlitz(t)$, $h=h(z)$, etc.:
$$g^\star_k=h^{q^k}(\tau^{k+1}s_\carlitz)\{\bsb{d}_1\tau^{k+1}(\bsb{d}_2)-\bsb{d}_2\tau^{k+1}(\bsb{d}_1)\}.$$
\end{Theorem}

Just as classical recurrent sequences, $\tau$-recurrent sequences
have {\em characteristic roots} (see Section \ref{recurrencesequences}). The interest of Theorem \ref{charroots} relies in that the characteristic roots of the $\tau$-recurrent sequence $(g^\star_k)_{k\geq 0}$
are explicitly computed, and turn out to be the functions $\bsb{d}_1,\bsb{d}_2.$

We will also compute series expansions of ``Eisenstein type", like (\ref{eq:definition_eisenstein}),
for the forms $g^\star_k$. To ease the next discussion, we mildly modify the aspect of the $g_k$'s when $k>0$:
$$
g_k(z)=(-1)^{k+1}\widetilde{\pi}^{1-q^k}[k]\cdots[1]\left(z\sideset{}{'}\sum_{c,d\in A}\frac{c}{(cz+d)^{q^k}}+\sideset{}{'}\sum_{c,d\in A}\frac{d}{(cz+d)^{q^k}}\right).
$$
The identity above can then be rewritten in a more compact form as a scalar product:
\begin{equation}\label{makeup}g_k(z)=-\zeta(q^k-1)^{-1}\mathcal{E}_k(z)\cdot\mathcal{F}_0(z),\end{equation}
where $\zeta$ is the restriction of Carlitz-Goss' zeta function at the integers, and where $\mathcal{E}_k(z)$ is the convergent (conditionally convergent 
when $k=0$)
series $\sum_{c,d\in A}'(c,d)(cz+d)^{-q^k}$, defining a holomorphic map $\Omega\rightarrow \mathbf{Mat}_{1\times 2}(\CC_\infty)$, and
$\mathcal{F}_0(z)$ is the map $\binom{z}{1}:\Omega\rightarrow \mathbf{Mat}_{2\times 1}(\CC_\infty)$.
It is not difficult to see that $\mathcal{E}_k,\mathcal{F}_0$ are {\em vectorial modular forms}. Similar obvious remarks can be 
made also for the classical Eisenstein series or Poincar\'e series, nonetheless, this will be one of the most important observations we make 
in this paper.

Let $t$ be an element of $\CC_\infty$. We have the ``evaluating at $t$" ring homomorphism
$$\chi_t:A\rightarrow \FF_q[t]$$
defined by $\chi_t(a)=a(t)$. In other words, $\chi_t(a)$ is the image of the polynomial map $a(t)$ obtained by substituting, in $a(\theta)$, $\theta$ by $t$. For example,
$\chi_t(1)=1$ and $\chi_t(\theta)=t$. The notation is motivated by the fact that if we choose $t\in\FF_q^{\text{alg.}}$ then $\chi_t$ factors through a Dirichlet character modulo the ideal generated by the minimal polynomial of $t$
in $A$.

Let $\alpha$ be a positive integer. We consider the integral value of $L$-series:
$$L(\chi_t,\alpha)=\sum_{a\in A^+}\chi_t(a)a^{-\alpha}=\prod_{\mathfrak{p}}(1-\chi_t(\mathfrak{p})\mathfrak{p}^{-\alpha})^{-1}\in K_\infty[[t]],$$
where $A^+$ denotes the set of monic polynomials of $A$,
converging for all $t\in B_q$ and not identically zero, and where the eulerian product runs over the monic irreducible 
polynomials of $A$.

\begin{Theorem}\label{scalarprod}
For all $z\in\Omega$, $t\in \CC_\infty$ with $|t|$ small enough, and $k\geq 0$, the following series expansion holds:
$$L(\chi_t,q^k)g_k^\star(z,t)=-\bsb{d}_1(z,t)\sideset{}{'}\sum_{c,d\in A}\frac{\chi_t(c)}{(cz+d)^{q^k}}-\bsb{d}_2(z,t)\sideset{}{'}\sum_{c,d\in A}\frac{\chi_t(d)}{(cz+d)^{q^k}}.$$
\end{Theorem}
Taking the limit $t\to\theta$ one recovers (\ref{makeup}) after some calculation.
For $t=\theta^{q^k}$, the series $L(\chi_t,q^k)$ diverges and the series on the right-hand side of the identity above
become conditionally convergent. However, it is possible to deduce from Theorem \ref{scalarprod} the following series expansion of ``Eisenstein type" for the modular forms $m_k$, in the style of (\ref{makeup}):
$$m_k=-\bsb{d}_1(z,\theta^{q^k})\sideset{}{'}\sum_{c,d\in A}\left(\frac{c}{cz+d}\right)^{q^k}-\bsb{d}_2(z,\theta^{q^k})\sideset{}{'}\sum_{c,d\in A}
\left(\frac{d}{cz+d}\right)^{q^k},$$
with a suitable choice of order of summation, because also these series are only conditionally convergent.

Theorem \ref{charroots} (and, up to a certain extent, Theorem \ref{scalarprod}) can be easily deduced by induction from 
the {\em deformation of Legendre's identity} (\ref{dethPsi}) that we copy here 
(\footnote{Proved in \cite{archiv}; the functions $g^\star_{k}$ already 
occur there.}):
$$h^{-1}(\tau s_\carlitz)^{-1}=\bsb{d}_1(\tau\bsb{d}_2)-\bsb{d}_2(\tau\bsb{d}_1).$$
In this paper however, we shall deduce Theorems \ref{charroots} and \ref{scalarprod} directly from 
a deeper result, Theorem \ref{firsttheorem} below, also implying (\ref{dethPsi}).

\medskip

\noindent\emph{Identities between deformations of vectorial modular forms.}
Both the series 
$$\bsb{e}_1(z,t)=\sideset{}{'}\sum_{c,d\in A}\frac{\chi_t(c)}{cz+d},\quad \bsb{e}_2(z,t)=\sideset{}{'}\sum_{c,d\in A}\frac{\chi_t(d)}{cz+d}$$
play a special role in Theorem \ref{scalarprod}. They
converge for $(z,t)\in\Omega\times B_q$ and define
functions $\Omega\rightarrow \CC_\infty[[t]]$ such that all the series in the images converge over $B_q$, and will be in the center of 
interest of this paper.

We will consider the vector function ${}^t\mathcal{E}$ with $\mathcal{E}=L(\chi_t,1)^{-1}(\bsb{e}_1,\bsb{e}_2)$ and $\mathcal{F}=\binom{\bsb{d}_1}{\bsb{d}_2}$ as examples of
{\em deformations of vectorial modular forms}. Thanks to this interpretation and the theory of $\tau$-linear recurrent sequences, we shall prove the following theorem, which immediately delivers Theorem \ref{charroots} (and Theorem \ref{scalarprod}):
\begin{Theorem}\label{firsttheorem}
The following identity holds in the domain $\Omega\times B_q$:
\begin{equation}\label{eqP101}
\mathcal{E}=(\tau s_\carlitz) h(-\tau\bsb{d}_2,\tau\bsb{d}_1).
\end{equation}
\end{Theorem}

\noindent\emph{Values of $L$-functions.} Theorem \ref{firsttheorem} implies the following corollary which, amazingly, does not seem to have been noticed previously:
\begin{Corollaire}\label{corollairezeta11}
The following identity holds:
$$L(\chi_t,1)=-\frac{\widetilde{\pi}}{\tau s_{\carlitz}}.$$ 
\end{Corollaire}

According to Corollary \ref{corollairezeta11}, the inverse of $\tau s_\carlitz=(t-\theta)s_\carlitz$ (that is, the function $\bsb{\Omega}$ of \cite{ABP}) is 
proportional to an $L$-value thus allowing entire analytic continuation in terms of the parameter $t$. 
It follows that
\begin{equation}\label{limit1}
\lim_{t\to\theta}L(\chi_t,1)=1\end{equation}
and the classical formulas for the values $\zeta(q^k-1)=\sum_{a\in A^+}a^{1-q^k}$ of Carlitz-Goss' zeta function can be easily deduced from this. Also, 
if $t=\xi\in\FF_q^\text{alg.}$, Corollary \ref{corollairezeta11} directly yields that $L(\chi_t,1)$, the value of an $L$-function associated to a Dirichlet character,
is a multiple of $\widetilde{\pi}$ by an algebraic element of $\CC_\infty$, solution of an algebraic equation $$X^{q^r-1}=(\xi-\theta^{q^{r-1}})\cdots(\xi-\theta).$$
Some of these consequences, but not all, are covered by the so-called Anderson {\em log-algebraic power series} identities for {\em twisted harmonic sums}, see \cite{Anbis,Anter,Lutes}, see also \cite{Dam} (\footnote{I am thankful to Vincent Bosser and to Matthew Papanikolas for having drawn my attention
to the papers \cite{Dam} and \cite{Lutes}, as well as for useful discussions concerning Corollary \ref{corollairezeta11}.}).

\medskip

\noindent\emph{The function $\bsb{E}$.}
Introduced in \cite{archiv}, this function is just the product: 
$$\bsb{E}=-h\tau\bsb{d}_2,$$ (see also \cite{BP3}), defined over $\Omega\times \CC_\infty$ and determining
a map $\Omega\rightarrow \CC_\infty[[t]]$ such that the series in the image have infinite radius of convergence; it has the 
property that $\bsb{E}(z,\theta)=E(z)$.
Theorem \ref{firsttheorem} implies that for $|t|<q$,\begin{equation}
\bsb{e}_1(z,t)=L(\chi_t,1)(\tau s_\carlitz)(t)\bsb{E}(z,t)=-\widetilde{\pi}\bsb{E}(z,t).\label{equ1}
\end{equation}
We then obtain the following result, providing a nice series expansion of $\bsb{E}$ ``near infinity":
\begin{Corollaire}\label{firstcorollary}
We have the following identity, valid for all $z\in\Omega$ and $t\in \CC_\infty$ such that $|t|<q^q$:
$$\bsb{E}(z,t)=\sum_{c\in A^+}\chi_t(c)u_c(z),$$
\end{Corollaire}
where we have used the functions $u_c(z):=e_\carlitz(c\widetilde{\pi}z)^{-1}$, with $e_\carlitz$ Carlitz's exponential.
For $t=\theta$, this reduces to \cite[(8.2)]{Ge}.

\medskip

\noindent\emph{Computation of $u$-expansions.}
Since for any $w$, every element $f$ of $M_{w,0}\otimes\FF_q[t]$ has $u$-expansion
$$f(z)=\sum_{i\geq 0}c_i(t)u(z)^i$$
with the $c_i$'s in $\CC_\infty[t]$, converging locally at $u=0$ for every fixed $t$, it is then very natural to
try to deduce from Theorem \ref{charroots} the $u$-expansions of the $g^\star_k$'s. Indeed, as we saw in \cite{archiv},
the function $\bsb{d}_2$ has a $u$-expansion with coefficients in $\FF_q[t,\theta]$, see (\ref{d2inu}). Unfortunately,
the function $s_\carlitz^{-1}\bsb{d}_1$ does not enjoy this property (see Lemma \ref{transcendence_of_d1}).

This paper contains a result, Theorem \ref{gkstarexpansion}, providing
a simple way to compute the $u$-expansions of the series $g^\star_k$'s
from the $u$-expansions of $\bsb{d}_2$ and a ``mysterious" function $\bsb{d}_3$ which allows $u$-expansion, introduced in Section 
\ref{computing}; the paper also presents 
simple algorithms to compute the $u$-expansion of $\bsb{d}_3$, as we did for $\bsb{d}_2$ in \cite{archiv}. We do not state 
Theorem \ref{gkstarexpansion} here (the statement requires some further preparation), but we mention a simple corollary of it.
\begin{Corollaire}\label{corgstar}
The truncation of the $u$-expansion of $g^\star_k$ to the order $q^k+2q-3$ is given, for $q\neq2$, by the truncation to the same order 
of the series:
$$\bsb{d}_2\left(1+\sum_{i=0}^{k-1}(t-\theta^{q^k})\cdots(t-\theta^{q^{i+1}})u^{q^k-q^i}\right)-(t-\theta)u^{q^k+q-2}.$$
\end{Corollaire}

The case $q=2$ is more involved, but can also be handled with the methods described here.
Choosing $t=\theta$ and using (\ref{limitsd}), we get the well known truncation to the order $q^k$ of the $u$-expansion of $g_k$ first computed by Gekeler in \cite{Ge}:
$$1+\sum_{i=0}^{k-1}(\theta-\theta^{q^k})\cdots(\theta-\theta^{q^{i+1}})u^{q^k-q^i}+\cdots.$$
 If on the other side we replace $t=\theta^{q^k}$ in the above expression, 
we obtain some coefficients of small order for the para-Eisenstein series $m_k$ provided we have knowledge 
of the $u$-expansion of $\bsb{d}_2$ up to a certain order, by using the algorithms developed in \cite{archiv}. 
Observe that then, the sum over $i=0,\ldots,k-1$ vanishes.

\medskip\noindent\emph{Link with extremal quasi-modular forms.}
To complete our paper, we will describe some links between the present work, $\cite{archiv}$, and the joint work \cite{BP2}.
For $l,w$ non-negative integers and $m$ a class of $\ZZ/(q-1)\ZZ$, we introduce the $\CC_\infty$-vector space of {\em Drinfeld quasi-modular forms
of weight $w$, type $m$ and depth $\leq l$}:
$$\widetilde{M}^{\leq l}_{w,m}=M_{w,m}\oplus M_{w-2,m-1}E\oplus\cdots\oplus M_{w-2l,m-l}E^l.$$

In \cite{BP2}, we have introduced the sequence of Drinfeld quasi-modular forms $(x_k)_{k\geq 0}$ with 
$x_k\in\widetilde{M}_{q^k+1,1}^{\leq 1}\setminus M$, defined by $x_0=-E$, $x_1=-Eg-h$ and by the recursion formula
\[x_k=x_{k-1}g^{q^{k-1}}-[k-1]x_{k-2}\Delta^{q^{k-2}},\quad k\geq 2,\] where we recall that $\Delta=-h^{q-1}$. The spaces $\widetilde{M}^{\leq l}_{w,m}$ embed in $\CC_\infty[[u]]$. In \cite[Theorem 1.2]{BP2}, we have showed that for all $k\geq 0$,
$x_k$ is {\em extremal}, in the sense that its order of vanishing at $u=0$, denoted by $\nu_\infty(x_k)$, is the biggest possible value for $\nu_\infty(f)$, if 
$f\in\widetilde{M}_{q^k+1,1}^{\leq 1}\setminus\{0\}$.
We also computed the order of vanishing: $\nu_\infty(x_k)=q^k$ for all $k$. 
After \cite[Proposition 2.3]{BP2}, the series expansion of 
$$E_k=(-1)^{k+1}\frac{x_k}{[1][2]\cdots[k]}$$   of $x_k$ begins with $u^{q^k}$
(where the empty product is $1$ by definition). Hence, with $E_0=E$, $E_k$ is the unique normalised extremal quasi-modular form
in $\widetilde{M}_{q^k+1,1}^{\leq 1}$ for all $k\geq 0$. 

We will obtain the following result.

\begin{Theorem}\label{propositionxk} For $k\geq 0$, we have 
\begin{eqnarray}\label{valueatthetaforE}E_k(z)&=&(\tau^k\boldsymbol{E})(z,\theta).
\label{firstET}
\end{eqnarray}
In particular, we have the series expansions
\begin{equation}\label{seriesEk}
E_k=\sum_{c\in A^+}cu_c^{q^k},
\end{equation} from which it is apparent that $E_k$ has $u$-expansion defined over
$A$.
\end{Theorem}

\noindent\emph{Remarks. 1.} The integrality of the coefficients of the normalised extremal quasi-modular form of weight $q^k+1$ and type $1$ 
supports Conjecture 2 of Kaneko and Koike in \cite{KK2}, asserting that if $f_{l,w}\in\QQ[[q]]$ is the $q$-expansion 
of the normalised extremal quasi-modular form of weight $w$ and depth $l\leq 4$, then $f_{l,w}\in\ZZ_p[[x]]$ for all $p>w$.

\medskip\noindent\emph{2.} 
The $u$-expansion (\ref{seriesEk}) follows from Corollary (\ref{firstcorollary}) applying (\ref{firstET}); hence, the integrality result 
is a consequence of our Theorem \ref{firsttheorem}.
It is easy to show that, for all $k\geq 0$, $g_k$ is the {\em extremal modular form} in $M_{q^k-1,0}$
(\footnote{This means that $g_k$ is the unique normalised form in $M_{q^k-1,0}\setminus\{0\}$ with the maximal order of vanishing of $g_k-1$.}). 

We also recall, from \cite{BP}, the derivation $D_1=u^2d/du$ on $\CC_\infty[[u]]$, 
which yields a $\CC_\infty$-linear map $\widetilde{M}^{\leq l}_{w,m}\rightarrow\widetilde{M}^{\leq l+1}_{w+2,m+1}$.
After Corollary \ref{corgstar}, the normalisation of $D_1g_k$ is the extremal quasi-modular form $E_k$:
\begin{equation}\label{secondET}
E_k=(-1)^{k}\frac{D_1g_k}{[1][2]\cdots[k]}.
\end{equation}
Roughly speaking,
in the classical theory of modular forms for the group $\mathbf{SL}_2(\ZZ)$ we have only one analogue 
of this striking situation, which is related to the {\em theta series} associated to the {\em Leech lattice} (of weight $12$):
$$\Theta_{\Lambda_{24}}=E_{12}-\frac{65520}{691}\Delta$$
where now, $E_{12}$ denotes the classical normalised classical Eisenstein series of weight $12$ and $\Delta$ is 
the normalised cusp form of weight $12$.

Let $f_{1,14}$ be the normalised extremal quasi-modular form of weight $14$ and depth $1$ in the sense of \cite{KK2}. 
The only analogue of the formula (\ref{secondET}) at $t=\theta$ in the classical framework is:
$$f_{1,14}=\frac{1}{393120}D\Theta_{\Lambda_{24}},$$ where $D$ denotes {\em Ramanujan's derivation} $(2\pi i)^{-1}d/dz$.
This agrees with the above mentioned conjecture of Kaneko and Koike because no prime exceeding $13$ divides $393120$.
 Numerical inspection suggests that $f_{1,14}$
is defined over $\ZZ$ but this property does not seem to be easy to prove
(\footnote{Notice however, that $f_{1,14}\in\ZZ_{11}[[q]]$ because $11$ does not divide $393120$. I am thankful to Gabriele Nebe 
for having observed that we also have $f_{1,14}\in\ZZ_{p}[[q]]$ for $p=5,7,13$ by using some properties of
the action of the double cover of the Conway group $2\text{Co}_1$ over $\Lambda_{24}$.}).
In the Drinfeldian case, the integrality of the coefficients of $E_k$ is an ultimate consequence 
of our formula (\ref{equ1}). So far, we do not know of an analogue of this formula in the classical framework.




\section{$\tau$-recurrent sequences\label{recurrencesequences}}

This section is devoted to the basic elements of the theory of $\tau$-recurrent sequences; the presentation is made in a
mild setting, yet more general than required by the rest of the paper.
In this section, $\mathcal{K}$ denotes any field endowed with an automorphism $\tau:\mathcal{K}\rightarrow \mathcal{K}$ of infinite order. We will refer to
the couple $(\mathcal{K},\tau)$ as to a {\em difference field}.
We denote by $\mathcal{K}^\tau$ the constant subfield of $\mathcal{K}$, that is, the subfield whose elements $x$ satisfy 
$\tau x=x$.

Let $x_1,\ldots,x_s$ be elements of $\mathcal{K}$. Their {\em $\tau$-wronskian} (sometimes called ``casoratian") is 
the determinant:
$$W_\tau(x_1,\ldots,x_s)=\det\left(\begin{array}{cccc} x_1 & \tau x_1 & \cdots & \tau^{s-1}x_1\\
x_2 & \tau x_2 & \cdots & \tau^{s-1}x_2\\
 \vdots & \vdots & & \vdots \\
 x_s & \tau x_s & \cdots & \tau^{s-1}x_s
\end{array}\right).$$
\begin{Lemme}\label{tauwronskian}
The elements $x_1,\ldots,x_s$ are $\mathcal{K}^\tau$-linearly independent if and only if $W_\tau(x_1,\ldots,x_s)\neq 0$.
\end{Lemme}
\noindent\emph{Proof.} This is a classical result that can be easily proved by induction on $s\geq 0$; we recall 
the proof here for convenience of the reader.
First of all, we notice that $W_\tau(x_1,\ldots,x_s)=0$ if and only if there exist elements $\lambda_1,\ldots,\lambda_s\in\mathcal{K}$,
not all zero, such that
\begin{equation}\label{lambdas}
x_1(\tau^n\lambda_1)+\cdots+x_s(\tau^n\lambda_s)=0,\quad n\in\ZZ.\end{equation}
Obviously, the lemma is true for $s=1$ so we consider $s>1$ and $x_1,\ldots,x_s$ such that $W_\tau(x_1,\ldots,x_s)=0$; there
exist $\lambda_1,\ldots,\lambda_s\in\mathcal{K}$ such that (\ref{lambdas}) holds.

If $\lambda_1,\ldots,\lambda_s$ are all in $\mathcal{K}^{\tau}$, we are done. Hence, we can assume that $\lambda_1\not\in\mathcal{K}^\tau$,
so that 
$$x_2(\tau^n\gamma_2)+\cdots+x_s(\tau^n\gamma_s)=0,\quad n\in\ZZ,$$
with $\gamma_i=\frac{\lambda_i}{\lambda_1}-\frac{\tau\lambda_i}{\tau\lambda_1}$ for $i=2,\ldots,s$ and the lemma 
follows by induction on $s$.\CVD

\noindent\emph{Remark.}
The proof is effective in the sense that the space generated by $(\lambda_1,\ldots,\lambda_s)$ can be explicitly 
computed in terms of the $x_i$'s following the inductive process step by step. Moreover, it sometimes (but not always) happens that limit processes, such as taking $n\rightarrow\infty$, and making use of some topology, furnish explicit
$\mathcal{K}^\tau$-linear dependence relations directly.

\medskip

We review now some elementary facts about {\em $\tau$-linear recurrent sequences} and their associated {\em $\tau$-linear equations}.

Let $L$ be a $\tau$-linear operator in the skew polynomial ring $\mathcal{K}[\tau]$. If $L=A_0\tau^0+\cdots+A_s\tau^s$ with $A_s\neq0$,
we will say that $L$ has order $s$. We will also say that the operator $L$ is {\em simple} if $A_0\neq0$. From now on, we will
only consider simple such operators, unless otherwise specified.

Let $\mathcal{G}:\ZZ\rightarrow \mathcal{K}$ be a sequence (this will be often denoted by $(\mathcal{G}_k)_{k\in\ZZ}$) and $\lambda\in \mathcal{K}$.
We shall write $\lambda * \mathcal{G}$ for the sequence
$$(\lambda * \mathcal{G})_{k\in\ZZ}=((\tau^k\lambda) \mathcal{G}_k)_{k\in\ZZ}.$$
With this action of $\mathcal{K}$, the set of the sequences $\ZZ\rightarrow \mathcal{K}$ is a vector space over $\mathcal{K}$.

Let $L=A_0\tau^0+\cdots+A_s\tau^s\in \mathcal{K}[\tau]$ be an operator as above, and let $\mathcal{G}$ be a sequence $\ZZ\rightarrow \mathcal{K}$. 
We 
will write $L(\mathcal{G})$ for the sequence 
$$L(\mathcal{G}):=(A_0\tau^0\mathcal{G}_k+\cdots+A_s\tau^s\mathcal{G}_{k-s})_{k\in\ZZ}.$$

We will say that $\mathcal{G}=(\mathcal{G}_k)_{k\in\ZZ}$ is a
{\em $\tau$-linear recurrent sequence} with coefficients in $\mathcal{K}$ associated to 
$L$ if 
\begin{equation}\label{taurecurrenceseq}
L(\mathcal{G})\equiv 0.\end{equation}
We will also say that $\mathcal{G}$ is of order $s$, if for any non-zero operator $L'\in \mathcal{K}[\tau]$ of order $<s$,
$L'(\mathcal{G})\not\equiv 0$.

Let $V=V(L)$ be the set of all the $\tau$-recurrent sequences satisfying (\ref{taurecurrenceseq}) for a given non-zero operator $L$ in $\mathcal{K}[\tau]$.
Since $L(\lambda * \mathcal{G})=\lambda * L(\mathcal{G})$, $V$ has a structure of $\mathcal{K}$-vector space; the dimension is finite, equal to $s$.

Assume that $V$ contains some constant sequence $(x)_{k\in\ZZ}$. Then, $x$ is a solution of  
the {\em associated linear $\tau$-difference equation} 
\begin{equation}\label{eqdiffinK}Lx=0.\end{equation}
In other words, with $L$ as above, we have $A_0x+A_1\tau x+\cdots +A_s \tau^s x=0$.
The set $V^\tau=V^\tau(L)$ of solutions of (\ref{eqdiffinK}) has a natural structure of $\mathcal{K}^\tau$-vector space.

\begin{Lemme}\label{r=s} Let $L$ be simple of order $s$ and let $V=V(L),V^\tau=V^\tau(L)$ be as above.
We have $\dim V^\tau\geq s$ if and only if $\dim V^\tau=s$. In the latter case, choose a basis $(x_1,\ldots,x_s)$ of $V^\tau$.
If $\mathcal{F}\in\mathbf{Mat}_{s\times 1}(\mathcal{K})$ is defined by ${}^t\mathcal{F}=(x_1,\ldots,x_s)$  (transpose),
then the map $$\mathcal{E}:V\rightarrow\mathbf{Mat}_{1\times s}(\mathcal{K})\cong \mathcal{K}^s$$ defined by
\begin{equation}\label{mathcalE}\mathcal{G}=(\mathcal{G}_k)_{k\in\ZZ}\mapsto \mathcal{E}(\mathcal{G}):=(\mathcal{G}_0,\tau^{-1}\mathcal{G}_1,\ldots,\tau^{-s+1}\mathcal{G}_{s-1})\cdot M^{-1}\end{equation}
is an isomorphism of $\mathcal{K}$-vector spaces.
\end{Lemme}

\noindent\emph{Proof.}
Let us assume that the dimension of $V^\tau$ is not smaller than $s$. Then, there exist $\mathcal{K}^\tau$-linear elements $x_1,\ldots,x_r$
of $\mathcal{K}$ solutions of (\ref{eqdiffinK}) with $r\geq s$.

Let $\mathcal{F}\in\mathbf{Mat}_{r\times 1}(\mathcal{K})$ be such that ${}^t\mathcal{F}=(x_1,\ldots,x_r)$.
By Lemma \ref{tauwronskian}, $W_\tau(x_1,\ldots,x_r)\neq 0$ and the matrix $M=(\mathcal{F},\tau^{-1}\mathcal{F},\ldots,\tau^{-r+1}\mathcal{F})$ is 
invertible.
The map \begin{equation}\label{mathcalE}V\rightarrow\mathbf{Mat}_{1\times r}(\mathcal{K})\cong \mathcal{K}^r
\end{equation} defined by
$$\mathcal{G}=(\mathcal{G}_k)_{k\in\ZZ}\mapsto (\mathcal{G}_0,\tau^{-1}\mathcal{G}_1,\ldots,\tau^{-r+1}\mathcal{G}_{r-1})\cdot M^{-1}$$
is then an isomorphism of $\mathcal{K}$-vector spaces and $r=s$. Therefore, $\dim_{\mathcal{K}^\tau}V^\tau=s$.\CVD


An operator $L=A_0\tau^0+\cdots+A_s\tau^s\in \mathcal{K}[\tau]$ of order $s$
is said to be {\em split} if $\dim V^\tau(L)=s$. A split operator is also simple. This definition obviously depends on the field $\mathcal{K}$.

\begin{Proposition}
Let $L\in \mathcal{K}[\tau]$ be a split operator of order $s$ and choose a basis $(x_1,\ldots,x_s)$ of $V^\tau(L)$. 
Let $V$ the $\mathcal{K}$-vector space of the $\tau$-recurrent sequences $\mathcal{G}$ such that $L(\mathcal{G})=0$.
Then, for all $\mathcal{G}\in V$ there exists one and only one element $\mathcal{E}\in\mathbf{Mat}_{1\times s}(\mathcal{K})$ such that
for all $k\in\ZZ$,\begin{equation}\label{EF}
\mathcal{G}_k=(\tau^k\mathcal{E})\cdot\mathcal{F}.
\end{equation}
\end{Proposition} 
\noindent\emph{Proof.} This follows from Lemma \ref{r=s}, taking $\mathcal{E}=\mathcal{E}(\mathcal{G})$ as in (\ref{mathcalE}).\CVD

\medskip

\begin{Proposition}\label{equationtau}
Let $x_1,\ldots,x_s$ be elements of $\mathcal{K}$.
Define, for $k=0,\ldots,s$, $$A_k=A_k^{\tau}(x_1,\ldots,x_s):=(-1)^{s+k}\det\left(\begin{array}{ccccc} \tau^0x_1 &  \cdots & \widehat{\tau^k x_1} & \cdots & \tau^{s}x_1\\
 \vdots & & \vdots & & \vdots \\
\tau^0x_s &  \cdots & \widehat{\tau^k x_s} & \cdots & \tau^{s}x_s
\end{array}\right),$$ where the hats mean that the corresponding column must be discarded.
Denote by $V^\tau(x_1,\ldots,x_s)$ the $\mathcal{K}^\tau$-vector space generated by the $x_i$ and let us consider the operator
\begin{equation}\label{theoperator}
L=L(x_1,\ldots,x_s)=A_0\tau^0+\cdots+A_s\tau^s.\end{equation}
Then,
$$V^\tau(L)=V^\tau(x_1,\ldots,x_s).$$
If the $x_i$'s are $\mathcal{K}^\tau$-linearly independent, then $L$ is split of order $s$.

Let $\mathcal{F}={}^t(x_1,\ldots,x_s)$ be a matrix of $\mathbf{Mat}_{s\times 1}(\mathcal{K})$ whose entries are $\mathcal{K}^\tau$-linearly independent. For all $\mathcal{E}\in\mathbf{Mat}_{1\times s}(\mathcal{K})$, the 
 sequence $\mathcal{G}=(\mathcal{G}_k)_{k\in\ZZ}$ defined by
 $$\mathcal{G}_k=(\tau^k\mathcal{E})\cdot\mathcal{F}$$ belongs to $V(L)$ with $L$ as in (\ref{theoperator}), 
 and every sequence of $V(L)$ can be expressed as above
 for some $\mathcal{E}$.
\end{Proposition}

\noindent\emph{Proof.} The existence of the operator $L$ follows easily by solving the $\tau$-difference equation
$$W_\tau(x_1,\ldots,x_s,X)=0.$$ Indeed,
by Lemma \ref{tauwronskian}, we have that $W_\tau(x_1,\ldots,x_s,x)=0$ with $x\in \mathcal{K}$ if and only if $x$ belongs to $\mathbf{Vect}_{\mathcal{K}^\tau}(x_1,\ldots,x_s)$. 
The non-vanishing of $A_0$ is also obvious as $A_0=(-1)^s\tau A_s=\tau W_\tau(x_1,\ldots,x_s)$.
 The final part of the proposition follows from a simple application of Proposition \ref{equationtau} which provides the 
operator $L$.\CVD
 
 The entries of $\mathcal{F}$ in Proposition \ref{equationtau} are called the {\em characteristic roots} of the $\tau$-linear recurrent sequence $\mathcal{G}_k$. 
 
 \medskip
 
Let us consider a sequence $\mathcal{G}$ as in (\ref{EF}), with $\mathcal{E},\mathcal{F}$ two matrices with 
entries in $\mathcal{K}$.
Then, with the above notations, we can introduce the {\em adjoint sequence} $\mathcal{H}=(\mathcal{H}_k)_{k\in\ZZ}$ defined by
$$\mathcal{H}_k=\tau^{-k}\mathcal{G}_k=\mathcal{E}\cdot(\tau^{-k}\mathcal{F}).$$

Let us assume that the entries of the matrix $\mathcal{E}\in\mathbf{Mat}_{1\times s}(\mathcal{K})$ are $\mathcal{K}^\tau$-linearly independent.
Then, $W_{\tau^{-1}}(\mathcal{E})\neq0$ and the above arguments with $\tau$ replaced by $\tau^{-1}$ ensure that $\mathcal{H}$ is a $\tau^{-1}$-recurrent sequence
of order $s$ so that there exists a split operator $L'\in \mathcal{K}[\tau^{-1}]$ of order $s$ such that $L'(\mathcal{H})=0$. If 
$L'=A_0'\tau^0+\cdots+A_s'\tau^s$ (so that $A_i'=A_i^{\tau^{-1}}(\mathcal{E}))$, then, for all $k\in \ZZ$
$$A_0'\mathcal{H}_k+A_1'\tau^{-1}\mathcal{H}_{k-1}+\cdots+A_s'\tau^{-s}\mathcal{H}_{k-s}=0.$$
Applying $\tau^s$ to the previous identities implies that the sequence $\mathcal{G}$ satisfies the following {\em $\tau$-linearised recurrent sequence} of order $s$:
$$
(\tau^{k}A_0')\mathcal{G}_k+(\tau^{k}A_1')\mathcal{G}_{k-1}+\cdots+(\tau^{k}A_s')\mathcal{G}_{k-s}=0,\quad k\in\ZZ.
$$
We will say that a sequence $\mathcal{G}$ of $\mathcal{K}$ is {\em generic} if there exist matrices $\mathcal{E}\in\mathbf{Mat}_{1\times s}(\mathcal{K})$
and $\mathcal{F}\in\mathbf{Mat}_{s\times 1}(\mathcal{K})$, both with $\mathcal{K}^\tau$-linearly independent entries, such that for all $k\in\ZZ$,
$$\mathcal{G}_k=(\tau^k\mathcal{E})\cdot\mathcal{F}.$$
Then, we obviously have the following proposition, containing all the properties encountered so far; later, we will use it for a specific generic
sequence of modular forms.
\begin{Proposition}\label{twotypessequence}
Let $\mathcal{G}=\mathcal{G}_k=(\tau^k\mathcal{E})\cdot\mathcal{F}$ be a generic sequence. If $L=A_0\tau^0+\cdots+A_s\tau^s$
is the split operator of $\mathcal{K}[\tau]$ associated to $\mathcal{F}$ and if $L'=A_0'\tau^0+\cdots+A_s'\tau^{-s}$
is the split operator of $\mathcal{K}[\tau^{-1}]$ associated to ${}^t\mathcal{E}$ (by Proposition
\ref{equationtau}), then $\mathcal{G}$ is at once $\tau$-linear recurrent and $\tau$-linearised recurrent (in both ways of order $s$). More precisely, for all $k\in\ZZ$,
\begin{eqnarray}
A_0\tau^0\mathcal{G}_k+A_1\tau^1\mathcal{G}_{k-1}+\cdots+A_s\tau^s\mathcal{G}_{k-s}&=&0,\\
(\tau^{k}A_0')\mathcal{G}_k+(\tau^{k}A_1')\mathcal{G}_{k-1}+\cdots+(\tau^{k}A_s')\mathcal{G}_{k-s}&=&0.
\end{eqnarray}
\end{Proposition} 

\subsection{Extending to existentially closed fields}\label{extending}

It is helpful, in some points of this paper, notably before computing 
solutions of certain linear $\tau$-difference equations, to first justify their existence in some extension field. 
In general, this makes no problem. However, we would like to point out here that by the so-called ``ACFA" theory of Chatzidakis and Hrushowski \cite{ChaHru}, 
there exist {\em existentially closed fields} $\KK$ containing $\mathcal{K}$ (more precisely, one speaks of the couple $(\KK,\tau)$
as being existentially closed). This means 
that there exists a field $\KK$ with an automorphism which extends $\tau$ (again denoted with $\tau$),
such that the constant subfield of $\KK$ for this automorphism is again $\mathcal{K}^\tau$, and such that every algebraic $\tau$-difference equation of positive order
has at least a non-zero solution $x\in \KK$.

\begin{Lemme} Let us assume that $(\mathcal{K},\tau)$ is existentially closed.
If $L=A_0\tau^0+\cdots+A_s\tau^s\in \mathcal{K}[\tau]$ is such that $A_sA_0\neq0$ as above,
then $\dim V^\tau(L)=s$.
\end{Lemme}
\noindent\emph{Proof.} 
We proceed by induction on $s\geq 0$. If $s=0$,
the statement of the lemma is trivial. Let us assume now that $s>0$.
Since $\mathcal{K}$ is existentially closed, there exists a solution $x_0\neq 0$ of $Lx=0$. Right division algorithm
holds in $\mathcal{K}[\tau]$, so that there exists $\tilde{L}\in \mathcal{K}[\tau]$ unique, with $L=\tilde{L}L_{x_0}$, where, for $x\in \mathcal{K}^\times$, we 
have written $L_x=\tau-(\tau x)/x$. Since the order of $\tilde{L}$ is $s-1$, there exist $y_1,\ldots,y_{s-1}$ $\mathcal{K}^\tau$-linearly independent elements of
$\mathcal{K}$ such that $\tilde{L}y_i=0$ for all $i$. Now, for all $i\geq 1$, let $x_i$ be a solution of $L_{x_0}x_i=y_i$
(they exist, again because $\mathcal{K}$ is existentially closed). Then, $x_0,x_1,\ldots,x_{s-1}$ are $s$ linearly independent 
elements of $\mathcal{K}$, solutions of $Lx=0$ so that $\dim V^\tau(L)\geq s$. By Lemma \ref{r=s}, $\dim V^\tau(L)=s$.\CVD

\noindent\emph{Remark.} The proof above requires that we solve non-homogeneous linear equations of order one as well, but this makes no problem even 
if were working in a field in which only linear homogeneous $n\times n$ systems of $\tau$-difference equations of order one can be solved, as we notice that 
the solution of the non-homogeneous equation $\tau x=ax+b$, reduces to the solution of the system $\tau y_1=ay_1+by_2$
and $\tau y_2=y_2$, which is homogeneous.
 
\section{Deformations of vectorial modular forms\label{dvmf}}

After having described some basic facts of the theory of $\tau$-linear recurrent sequences, we come back to our modular forms and we now 
start dealing with vectorial modular forms and their deformations. For this, we are making again specific choices of $\mathcal{K},\tau$ etc.

\subsection{Notation, tools}

Let $t$ be an indeterminate indipendent on $\theta$.
Often in this paper, $t$ will be also a parameter varying in $\CC_\infty$ and we will freely switch from formal series to functions. 

For a positive real number $r$, we denote by $\TT_{<r}$ the sub-$\CC_\infty$-algebra of $\CC_\infty[[t]]$ whose elements are formal series $\sum_{i\geq 0}c_it^i$ that converge 
for any $t\in \CC_\infty$ with $|t|<r$. 
We also denote by $\TT_\infty$ the sub-$\CC_\infty$-algebra of series that converge everywhere in $\CC_\infty$.
If $r_1>r_2>0$, we have $$\TT_{<r_2}\supset\TT_{< r_1}\supset\TT_\infty.$$ The {\em Tate algebra} 
of formal series of $\CC_\infty[[t]]$ converging for all $t$ such that $|t|\leq 1$ will be denoted by $\TT_1$ or $\TT$; it is contained in $\TT_{<1}$ and contains
$\TT_{<1+\epsilon}$ for all $\epsilon>0$; clearly, $\CC_\infty[[t]]\supset\TT_1\supset\TT_\infty$.

The ring $\CC_\infty[[t]]$ is  endowed with the $\FF_q[[t]]$-linear automorphism $\tau$ acting on formal series as follows:
$$\tau\sum_ic_it^i=\sum_ic_i^qt^i.$$ This automorphism induces automorphisms of $\TT_1,\TT_\infty$.

We will work with certain functions $f:\Omega\times B_r\rightarrow \CC_\infty$ with 
the property that for all $z\in\Omega$,
$f(z,t)$ can be identified with and element of $\TT_{<r}$. For such functions we will then also write
$f(z)$ to stress the dependence on $z\in\Omega$ when we want to consider them as functions $\Omega\rightarrow\TT_{<r}$ for some $r$.
Sometimes, we will not specify the variables $z,t$ and just write $f$ instead of $f(z,t)$ or $f(z)$ to lighten our formulas. Moreover,
$z$ will always denote a variable in $\Omega$ of any modular form in this paper.

In all the following, $\mathbf{Hol}(\Omega)$ denotes the ring of holomorphic functions on $\Omega$ and
$\mathbf{Me}(\Omega)$ its fraction field.
For $r$ a positive real number, let us denote by 
$\mathcal{R}_{<r}$ (resp. $\mathcal{R}$ or $\mathcal{R}_1$) the (integral) ring whose elements are the formal series $f=\sum_{i\geq0}f_it^i$,
such that
\begin{enumerate}
\item For all $i$, $f_i$ is a map $\Omega\rightarrow \CC_\infty$ belonging to $\mathbf{Hol}(\Omega)$.
\item For all $z\in\Omega$, $\sum_{i\geq0}f_i(z)t^i$ is an element of $\TT_{<r}$ (resp. $\TT$).
\end{enumerate}
We shall write $$\mathcal{R}_\infty=\bigcap_{r>0}\mathcal{R}_{<r}$$ and allow $r$ to vary in $\RR_{>0}\cup\{\infty\}$. The fraction fields $\mathcal{L}$ and $\mathcal{L}_\infty$ of the rings $\mathcal{R}$ and $\mathcal{R}_\infty$
are endowed with injective endomorphisms $\tau$ acting on formal series as follows:
$$\tau\sum_{i\geq0}f_i(z)t^i=\sum_{i\geq0}f_i(z)^qt^i.$$
\begin{Lemme}\label{fresnelvanderput}
We have $\mathcal{L}^\tau=\mathcal{L}_\infty^\tau=\FF_q(t)$.
\end{Lemme}
\noindent\emph{Proof.} It suffices to compute $\mathcal{L}^\tau$. Let $f=\sum_{i}f_i(z)t^i$ be in $\mathcal{L}^\tau$ and let us choose $z\in \Omega$
such that $\phi_i=f_i(z)$ is well defined for all $i$. We get a series $\phi_z=\sum_{i}\phi_it^i$ in $\LL^\tau$, where $\LL$ is the fraction 
field of $\TT$. It is well known (see \cite{Pa} or use use Theorem 2.2.9 of \cite{FP}) that $\LL^\tau=\FF_q(t)$.
This means that for all $z\in\Omega$ such that $f(z)$ is well defined, $f(z)=\phi_z\in\FF_q(t)$. Since the functions $f_i$ are 
meromorphic, we then get $f\in\FF_q(t)$.\CVD

After Section \ref{extending}, the fields $\mathcal{L}$ and $\mathcal{L}_\infty$ can be embedded in a field $\mathcal{K}$ endowed
with an automorphism extending $\tau$ (and denoted again by $\tau$), such that $(\mathcal{K},\tau)$ is existentially closed,
with $\mathcal{K}^\tau=\FF_q(t)$. From now on, we will apply the results of Section \ref{recurrencesequences}
to this difference field.

We end this preparatory section with some conventions on $u$-expansions.
We will say that a series $\sum_{i\geq i_0}c_iu^i$ (with the coefficients $c_i$ in some ring) is {\em normalised}, if $c_{i_0}=1$.
We will also say that the series is {\em of type} $m\in\ZZ/(q-1)\ZZ$ if $i\not\equiv m\pmod{q-1}$ implies $c_i=0$.
This definition is obviously compatible with the notion of type of a Drinfeld modular form already discussed in the introduction.

\subsection{Basic properties of vectorial modular forms.}

In this subsection we introduce {\em deformations of vectorial modular forms}. 
This part is largely inspired by a conspicuous collection of papers about vectorial modular forms for $\mathbf{SL}_2(\ZZ)$ notably 
by Knopp, Mason.

To make our list of references self contained, we only mention \cite{KM, Mas}, leaving the reader to further explore the literature.
In particular, we learned from \cite[Section 3]{KM} how to construct {\em vectorial Poincar\'e series}, of which we propose a Drinfeldian
counterpart in Subsection \ref{poincare}. It should be noticed, however, that our construction is not a complete adaptation
of Knopp and Mason's constructions and the analogy is superficial. Indeed, there is a fundamental gap between the theories. 
While symmetric powers of two-dimensional irreducible representations of 
$\SL_2(\ZZ)$ are irreducible, symmetric powers of two-dimensional irreducible representations of 
$\GL_2(A)$ are the most often {\em not} irreducible and split along tiny irreducible sub-representations. 

\medskip

Let us consider a representation
\begin{equation}\label{rho}
\rho:\Gamma\rightarrow \mathbf{GL}_s(\FF_q((t))).
\end{equation}
We assume that the determinant representation $\det(\rho)$ is the 
$\mu$-th power of the determinant character, for some $\mu\in\ZZ/(q-1)\ZZ$. In all the following, given $\gamma\in\Gamma$,
we denote by $J_\gamma$ the associated factor of automorphy $(\gamma,z)\mapsto cz+d$, if $\gamma=\sqm{a}{b}{c}{d}$.

\begin{Definition}
{\em A {\em deformation of vectorial modular form} (abridged to DVMF) of weight $w$, dimension $s$, type $m$ and radius $r\in\RR_{>0}\cup\{\infty\}$ associated with a representation $\rho$ as in 
(\ref{rho}) is a column matrix $\mathcal{F}\in\mathbf{Mat}_{s\times 1}(\mathcal{R}_{<r})$ such that, considering $\mathcal{F}$ as a map $\Omega\rightarrow\mathbf{Mat}_{s\times 1}(\TT_r)$ we have,
for all $\gamma\in\Gamma$, 
$$\mathcal{F}(\gamma(z))=J_\gamma^w\det(\gamma)^{-m}\rho(\gamma)\cdot\mathcal{F}(z).$$
The definition means that if the radius is $\infty$, then the entries of $\mathcal{F}$ are in $\mathcal{R}_\infty$.}
\end{Definition}
The set of deformations of vectorial modular forms of weight $w$, dimension $s$, type $m$ and radius $r$ associated to a representation $\rho$
is a $\TT_{<r}$-module (or $\TT_\infty$-module if $r=\infty$) that we will denote by $\mathcal{M}^{s}_{w,m}(\rho,r)$ or $\mathcal{M}^{s}_{w,m}(\rho)$ when the 
reference to a particular radius is clear. 

In this paper, $M^!_{w,m}$ denotes the $\CC_\infty$-vector space (of infinite dimension)
generated by quotients $f/g$ with $f\in M_{w',m'}$, $g\in M_{w'',m''}\setminus\{0\}$ such that $w'-w''=w$ and $m'-m''=m$.

If $s=1$ and if $\rho=\bsb{1}$ is the constant map, then $\mathcal{M}^1_{w,m}(\bsb{1},r)=M^!_{w,m}\otimes\TT_{<r}$. Therefore, for general $s$,
we have
a graded $M^!_{w,m}\otimes\TT_{<r}$-module $$\mathcal{M}^s(\rho,r)=\bigoplus_{w,m}\mathcal{M}^{s}_{w,m}(\rho,r).$$

\begin{Lemme}\label{twist} Let $k$ be a non-negative integer.
If $\mathcal{F}$ is in $\mathcal{M}^{s}_{w,m}(\rho,r)$, then $\tau^k\mathcal{F}\in\mathcal{M}^{s}_{wq^k,m}(\rho,r^{q^k})$
and $(\tau^{-k}\mathcal{F})^{q^k}\in\mathcal{M}^{s}_{w,m}(\rho,r)$.
\end{Lemme}

\noindent\emph{Proof.} from the definition, 
$$(\tau^k\mathcal{F})(\gamma(z))=J_\gamma^{wq^k}\det(\gamma)^{-m}\rho(\gamma)(\tau^k\mathcal{F})$$
because $\tau(\rho(\gamma))=\rho(\gamma)$.\CVD

\begin{Proposition}\label{wronsk} Let us assume that $r>1$, let 
us consider $\mathcal{F}$ in $\mathcal{M}^{s}_{w,m}(\rho,r)$
and 
${}^t\mathcal{E}$ in $\mathcal{M}^{s}_{w',m'}({}^t\rho^{-1},r)$,
choose nonnegative integers $k_1,\ldots,k_s$ and set 
$k=\max\{k_1,\ldots,k_s\}$. Then
$$\det(\tau^{k_1}\mathcal{F},\ldots,\tau^{k_s}\mathcal{F})\in M^!_{w(q^{k_1}+\cdots+q^{k_s}),sm+\mu}\otimes\TT_{<r},$$
and
$$\det(\tau^{-k_1}({}^t\mathcal{E}),\ldots,\tau^{-k_s}({}^t\mathcal{E}))^{q^k}\in M^!_{w'(q^{k-k_1}+\cdots+q^{k-k_s}),sm'-\mu}\otimes\TT_{<r}.$$
In particular, $$W_\tau(\mathcal{F})\in M^!_{w(1+q+q^2+\cdots+q^{s-1})),sm+\mu}\otimes\TT_{<r},$$
and $$W_{\tau^{-1}}({}^t\mathcal{E})^{q^{s-1}}\in M^!_{w'(1+q+q^2+\cdots+q^{s-1})),sm'-\mu}\otimes\TT_{<r}.$$
Moreover, for nonnegative $k$, if $\mathcal{G}_k$ denotes $(\tau^k\mathcal{E})\cdot\mathcal{F}$, then
$$\mathcal{G}_k\in M^!_{w+w'q^k,m+m'}\otimes\TT_{<r}.$$ 
\end{Proposition}
\noindent\emph{Proof.} Define the matrix function:
$$\mathbf{M}_{k_1,\ldots,k_s}=(\tau^{k_1}\mathcal{F},\ldots,\tau^{k_s}\mathcal{F}).$$
After Lemma \ref{twist} we have, for $\gamma\in\mathbf{GL}_2(A)$:
$$\mathbf{M}_{k_1,\ldots,k_s}(\gamma(z))=\det(\gamma)^{-m}\rho(\gamma)\cdot \mathbf{M}_{k_1,\ldots,k_s}(z)\cdot\mathbf{Diag}(J_\gamma^{wq^{k_1}},\cdots,
J_\gamma^{wq^{k_s}}).$$

If the $k_i$'s are all positive, the coefficients of the $t$-expansions of the entries of $\det(\tau^{k_1}\mathcal{F},\ldots,\tau^{k_s}\mathcal{F})$ are holomorphic functions on $\Omega$.
It they are all negative, the corresponding coefficients, raised to the power $q^k$, are holomorphic on $\Omega$.

The part of the proposition involving the determinant of $\mathbf{M}_{k_1,\ldots,k_s}$ follows easily. There is no additional 
difficulty in proving the part concerning the form $\mathcal{E}$. 

Also the latter property of the sequence $(\mathcal{G}_k)_k$
follows easily from Lemma  \ref{twist}. Indeed,
by this lemma, $\tau^k({}^t\mathcal{E})$ is in $\mathcal{M}^s_{wq^k,m'}({}^t\rho^{-1},r)$. Let 
$\gamma$ be in $\mathbf{GL}_2(A)$. We know that
$$(\tau^k\mathcal{E})(\gamma(z))=J_\gamma^{wq^k}\det(\gamma)^{-m}\;{}^t\mathcal{E}(z)\cdot\rho^{-1}(\gamma)$$
and
$$\mathcal{F}(\gamma(z))=J_\gamma^{w'}\det(\gamma)^{-m'}\rho(\gamma)\cdot\mathcal{F}(z).$$
Hence,
$$\mathcal{G}_k(\gamma(z))=J_\gamma^{wq^k+w'}\det(\gamma)^{-m-m'}\mathcal{G}_k(z),$$
from which we deduce that $\mathcal{G}_k\in M_{wq^k+w',m+m'}^!\otimes \TT_{<r}$.\CVD

The next proposition is a mere reproduction of the main properties described in Section \ref{recurrencesequences}
in the framework of deformations of vectorial modular forms.

\begin{Proposition}\label{equationf}
Assuming that $r>1$, let us consider $\mathcal{F}$ in $\mathcal{M}^{s}_{w,m}(\rho,r)$ and
let $\mathcal{E}$ be such that ${}^t\mathcal{E}$ is in $\mathcal{M}^{s}_{w',m'}({}^t\rho^{-1},r)$.
For $k\in\ZZ$, let us write $\mathcal{G}_k=(\tau^k\mathcal{E})\cdot\mathcal{F}$.

Then, for all $k=0,\ldots,s$, we have $$A_k=A_k^\tau(\mathcal{F})\in M^!_{(1+q+\cdots+\widehat{q^k}+\cdots+q^{s})w,sm+\mu}\otimes\TT_{<r}$$
(the hat means that we skip the corresponding term in the sum).
Let $L$ be the operator $A_0\tau^0+\cdots+A_s\tau^{s}$. If $r>1$ and if the components of $\mathcal{F}$ are $\FF_q(t)$-linearly independent, then 
$L$ is split, for any $k$ nonnegative integer,
$\mathcal{G}_k$ is an element of $M_{wq^k+w',m+m'}^!\otimes \TT_{<r}$
and we have the relations
\begin{eqnarray*}
A_0\mathcal{G}_k+A_1\tau \mathcal{G}_{k-1}+\cdots+A_s\tau^s\mathcal{G}_{k-s}&=&0.
\end{eqnarray*}

For all $k=0,\ldots,s$, $A_k'=A_k^{\tau^{-1}}(\mathcal{E})$
is such that $$(A_k')^{q^s}\in M^!_{(1+q+\cdots+\widehat{q^k}+\cdots+q^{s})w,sm-\mu}\otimes\TT_{<r}.$$ If the entries of $\mathcal{E}$ are
$\FF_q(t)$-linearly independent and $r>1$, then the operator $L'=A_0'\tau^0+\cdots+A_s'\tau^{-s}$ is split and
we also have the relations:
\begin{eqnarray*}
(\tau^kA_0')\mathcal{G}_k+(\tau^kA_1')\mathcal{G}_{k-1}+\cdots+(\tau^kA_s')\mathcal{G}_{k-s}&=&0.
\end{eqnarray*}
\end{Proposition}
\noindent\emph{Proof.} By Lemma \ref{fresnelvanderput} and Lemma \ref{tauwronskian}, the $\tau$-wronskian of 
$\mathcal{F}$ is non-zero. We apply Proposition \ref{equationtau} to obtain 
that $L$ is split of order $s$ and if the components of $\mathcal{E}$ are $\FF_q(t)$-linearly
independent, also $L'$ is split. By Proposition \ref{wronsk}, the coefficients $A_i$ are 
modular as claimed. The part of the proposition involving properties of the form $\mathcal{E}$ is similar and left to the reader.
Then, the proof of the proposition can be completed with the help of
Proposition \ref{twotypessequence}.\CVD

\subsection{Examples.}

From now on, we will use the representation $\rho=\rho_{t,1}:\mathbf{GL}_2(A)\rightarrow \mathbf{GL}_2(\FF_q[t])$ defined 
by 
$$\rho_{t,1}(\gamma)=\sqm{\ol{a}}{\ol{b}}{\ol{c}}{\ol{d}}$$ if $\gamma=\sqm{a}{b}{c}{d}\in\Gamma$, 
and its {\em symmetric powers} of order $l$ for $l\geq 1$
$$\rho_{t,l}=S^l(\rho_{t,1}):\mathbf{GL}_2(A)\rightarrow \mathbf{GL}_{l+1}(\FF_q[t]),$$
realised in the space of polynomial homogeneous of degree $s=l+1$ with coefficients in $\FF_q[t]$:
$$\rho_{t,l}\left(\sqm{a}{b}{c}{d}\right)(X^{s-r}Y^r)=(\ol{a}X+\ol{c}Y)^{s-r}(\ol{b}X+\ol{d}Y)^r.$$
The determinant
of $\rho_{t,l}$ is the $l(l+1)/2$-th power of the determinant character:
$$\det(\rho_{t,l}(\gamma))=\det(\gamma)^{\frac{l(l+1)}{2}}.$$
Together with $\rho_{t,l}$ we will also use the representation ${}^t\rho_{t,l}^{-1}$ (transpose of the inverse)
and we set $\rho_{t,0}(\gamma)=1$ for all $\gamma$. 

\subsubsection{First example: the functions $\Phi_l$}
We first discuss again the functions $\bsb{d}_1,\bsb{d}_2$ mentioned in the introduction.

For $z\in\Omega$, we have denoted by $\Lambda_z$ the $A$-module $A+zA$,
and we have the expression (\ref{ellipticexponential}) for the exponential function $e_{\Lambda_z}$.
We recall that:
\begin{eqnarray*}
\bfs_1(z,t)&=&\sum_{i=0}^\infty\frac{\alpha_i(z)z^{q^i}}{\theta^{q^i}-t}\\
\bfs_2(z,t)&=&\sum_{i=0}^\infty\frac{\alpha_i(z)}{\theta^{q^i}-t}.
\end{eqnarray*} These are
functions $\Omega\times B_q\rightarrow \CC_\infty$. From \cite{archiv}, we deduce that $\bsb{s}_1,\bsb{s}_2$ lie in $\mathcal{R}_{<q}$.

At $\theta$, the functions $\bfs_i(z,\cdot)$ have simple poles. Their respective residues are $-z$ for the function $\bfs_1(z,\cdot)$ and $-1$
for $\bfs_2(z,\cdot)$.
Moreover, we have $\bfs_1^{(1)}(z,\theta)=\eta_1$ and $\bfs_2^{(1)}(z,\theta)=\eta_2$, where $\eta_1,\eta_2$
are the {\em quasi-periods} of $\Lambda_z$ (see \cite[Section 4.2.4]{Bourbaki} and \cite[Section 7]{gekeler:compositio}).
We set, for $i=1,2$:
$$\bsb{d}_i(z,t):=\widetilde{\pi}s_\carlitz(t)^{-1}\bfs_i(z,t),$$ with $s_\carlitz$ defined in (\ref{scarlitz}).
We point out that, in the notations of \cite{archiv}, $\bsb{d}=\bsb{d}_2$. At first sight, we only have $\bsb{d}_1,\bsb{d}_2\in\mathcal{R}_{<q}$.
However, one sees easily that $s_\carlitz^{-1}\in\TT_\infty$ from which it follows that $\bsb{d}_1,\bsb{d}_2\in\mathcal{R}_\infty$.

The functions $\bsb{d}_1,\bsb{d}_2$ enjoy several properties that can be easily deduced from \cite{archiv}. Here, we are concerned with
a $\tau$-difference linear equation, a {\em deformation of Legendre's identity}, the quality of being 
a deformation of vectorial modular form and a $u$-expansion for $\bsb{d}_2$. 
These properties 
where obtained in \cite{archiv} for the functions $\bsb{s}_1,\bsb{s}_2$. Here we collect them in the following 
proposition, in terms of the functions $\bsb{d}_1,\bsb{d}_2$ (the deduction of the proposition from \cite{archiv} is immediate).

\begin{Proposition}\label{prarchiv}
We have five properties for the $\bsb{d}_i$'s.

\begin{enumerate}
\item $\bsb{d}_1,\bsb{d}_2\in\mathcal{R}_{\infty}$.
\item Let us write $\Phi_1=\binom{\bsb{d}_1}{\bsb{d}_2}$. We have 
$\Phi_1\in\mathcal{M}^2_{-1,0}(\rho_{t,1},\infty)$.
\item The following $\tau$-linear difference equations hold:
\begin{equation}\label{equexpansion}
\bsb{d}_i=(t-\theta^q)\Delta\bsb{d}_i^{(2)}+g\bsb{d}_i^{(1)},\quad i=1,2.
\end{equation}
\item Let us consider the matrix function:
$$\Psi (z,t):=\sqm{\bsb{d}_1(z,t)}{\bsb{d}_2(z,t)}{\bsb{d}_1^{(1)}(z,t)}{\bsb{d}_2^{(1)}(z,t)}.$$
For all $z\in\Omega$ and $t$ with $|t|<q$:
\begin{equation}\label{dethPsi}
\det(\Psi)=(t-\theta)^{-1}h(z)^{-1}s_\carlitz(t)^{-1}.
\end{equation}
\item We have the series expansion
\begin{equation}\label{uexpd2}
\bsb{d}_2 =\sum_{i\geq 0}c_i(t)u^{(q-1)i}\in1+u^{q-1}\FF_q[t,\theta][[u^{q-1}]],
\end{equation}
convergent for $t,u$ sufficiently close to $(0,0)$.
\end{enumerate}
\end{Proposition}

More precisely, we showed in \cite{archiv} that the series expansion in powers of $u$ of $\bsb{d}_2$ is as follows:
\begin{equation}\label{d2inu}
\bsb{d}_2=1+(\theta-t)u^{q-1}+(\theta-t)u^{(q^2-q+1)(q-1)}+\cdots\in1+(t-\theta)\FF_q[t,\theta][[u^{q-1}]],\end{equation}
where the dots $\cdots$ stand for terms of higher degree in $u$. 

For $l\geq 1$ fixed, let us consider the function:

$$\Phi_l=\left(\begin{array}{c}\bsb{d}_1^l\\ \bsb{d}_1^{l-1}\bsb{d}_2\\ \vdots\\ \bsb{d}_1\bsb{d}_2^{l-1}\\ 
\bsb{d}_2^l\end{array}\right):\Omega\rightarrow\mathbf{Mat}_{l+1\times 1}(\TT_\infty),$$ so that $\Phi_l\in\mathbf{Mat}_{l+1\times 1}(\mathcal{R}_\infty)$.

\begin{Lemme}\label{componentslinearlyindep}
We have $\Phi_{l}\in\mathcal{M}^{l+1}_{-l,0}(\rho_{t,l},\infty)$ and
the components of $\Phi_l$
are $\FF_q(t)$-linearly independent.
\end{Lemme}
\noindent\emph{Proof.} The first property is obvious after Proposition \ref{prarchiv}. Assume that we have a non-trivial linear dependence relation
with the $c_i$'s in $\FF_q(t)$:
$$\sum_{i=0}^{l}c_i\bsb{d}_1^i\bsb{d}_2^{l-i}=0.$$
Then, replacing $t=\theta$ and using (\ref{limitsd}), 
we find
$$\sum_{i=0}^lc_i(\theta)z^i=0$$
which is impossible.\CVD

By Proposition \ref{wronsk} and Lemma \ref{componentslinearlyindep},
there is a split operator of order $l+1$:
\begin{equation}\label{eqoperator}
L_l=A_{l,0}\tau^0+\cdots+A_{l,l+1}\tau^{l+1}\end{equation}
such that $L_l\Phi_l=0.$
In particular, $A_{l,l+1}\neq 0$. Moreover, it is 
easy to check that for all $l$ and $0\leq i\leq l+1$, there exists an integer $\mu=\mu(i,l)$ such that
$h^\mu A_{l,i}\in M_{*,*}\otimes \TT$ (use (\ref{dethPsi})).

\medskip

\noindent\emph{Examples.}
More specifically, if $l=1$, we find
$$A_{1,2}=\det(\Psi),\quad A_{1,1}=-\frac{\det(\Psi)g}{\Delta(t-\theta^q)},\quad A_{1,0}=-\frac{\det(\Psi)}{\Delta(t-\theta^q)},$$ 
implying (\ref{equexpansion}), and 
\begin{equation}\label{ell1}
L_1=-\tau^0-g\tau+\Delta(t-\theta^q)\tau^2.\end{equation}

If $l=2$, we find, after some rather heavy computation using (\ref{equexpansion}) and (\ref{dethPsi}):

\begin{eqnarray*}
W_\tau(\Phi_2)=A_{2,3}&=&\frac{\det(\Psi)^3g}{\Delta^2(\theta^q-t)^2},\\
A_{2,2}&=&\frac{\det(\Psi)^3g^q(g^{1 + q} + \Delta (t - \theta^q))}{\Delta^{
          2 + 2q}(\theta^q - t)^2(\theta^{q^2} - t)^2}\\
A_{2,1}&=&-\frac{\det(\Psi)^3g(g^{1 + q} + \Delta (t - \theta^q))}{\Delta^{
          3 + 2q}(\theta^q - t)^3(\theta^{q^2} - t)^2},\\
A_{2,0}&=&-\frac{\det(\Psi)^3g^q}{\Delta^{3+2q}(\theta^q-t)^3(\theta^{q^2}-t)^2},
\end{eqnarray*}
and 
\begin{equation}\label{ell2}
L_2=-\tau^0-g^{1-q}(g^{1+q}+\Delta(t-\theta^q))\tau+(g^{1+q}+\Delta(t-\theta^q))\Delta(\theta^q-t)\tau^2+g^{1-q}\Delta^{1+2q}(\theta^q-t)(\theta^{q^2}-t)^2\tau^3.\end{equation}

The explicit determination of the coefficients of the operator (\ref{eqoperator})
for the vectorial forms $\Phi_l$ for general $l$ looks like a difficult computational problem.

\subsubsection{Second example: Deformations of vectorial Poincar\'e series\label{poincare}}

Following \cite{Ge}, let us consider the subgroup $H=\left\{\sqm{*}{*}{0}{1}\right\}$ of $\Gamma=\mathbf{GL}_{2}(A)$ and its left action on $\Gamma$.

For $\delta=\sqm{a}{b}{c}{d}\in\Gamma$, the map $\delta\mapsto(c,d)$ induces a bijection
between the orbit set $H\backslash\Gamma$ and the set of $(c,d)\in A^2$ with $c,d$ relatively prime. 
For $l\geq 0$, let $V_l (\delta)$ be the row matrix $$(\ol{c}^l,\ol{c}^{l-1}\ol{d},\ldots,\ol{c}\ol{d}^{l-1},
\ol{d}^l).$$

We consider the factor of automorphy $$\mu_{\alpha,m}(\delta,z)=\det(\delta)^{-m}J_\gamma^\alpha,$$ where $m$ and $\alpha$ are 
positive integers (later, $m$ will also determine a type, that is, a class modulo $q-1$).

It is easy to show that the quantity $$\mu_{\alpha,m}(\delta,z)^{-1}u^m(\delta(z))V_l (\delta)$$ only depends on the class of $\delta\in H\backslash\Gamma$,
so that we can consider the following series:
$$\mathcal{E}_{\alpha,m,l}(z)=\sum_{\delta\in H\backslash\Gamma}\mu_{\alpha,m}(\delta,z)^{-1}u^m(\delta(z))V_l (\delta),$$ 
which is a row matrix  whose  $l+1$ entries are formal series. 

Let $\mathcal{V}$ be the set of functions $\Omega\rightarrow\mathbf{Mat}_{1\times l+1}(\CC_\infty[[t]]).$
We introduce, for $\alpha,m$ integers, $f\in \mathcal{V}$ and $\gamma\in\Gamma$, the Petersson slash operator:
$$f|_{\alpha,m}\gamma=\det(\gamma)^{m}(cz+d)^{-\alpha}f(\gamma(z))\cdot\rho_{t,l}(\gamma).$$
This will be used in the next proposition, where $\log^+_q(x)$ denotes the maximum between $0$ and $\log_q(x)$, the logarithm in base $q$  of $x>0$.

\begin{Proposition}\label{mainproppoincare} Let $\alpha,m,l$ be non-negative integers with $\alpha\geq 2m+l$, and write $r(\alpha,m,l)=(\alpha-2m-l)/l$
if $l\neq0$, and $r(\alpha,m,l)=\infty$ if $l=0$. We have the following properties.

\begin{enumerate}

\item For $\gamma\in\Gamma$, the map $f\mapsto f|_{\alpha,m}\gamma$ induces a permutation of the subset of $\mathcal{V}$: 
$$\mathcal{S}=\{\mu_{\alpha,m}(\delta,z)^{-1}u^m(\delta(z))V_l (\delta);\delta\in H\backslash\Gamma\}.$$

\item If $t\in \CC_\infty$ is chosen so that $r(\alpha,m,l)>\log^+_q|t|$, then the components of $\mathcal{E}_{\alpha,m,l}(z,t)$ are series of functions
of $z\in\Omega$ which converge absolutely and uniformly on every compact subset of $\Omega$ to holomorphic functions.

\item If $\log_q|t|<0$, then the components of $\mathcal{E}_{\alpha,m,l}(z,t)$ converge absolutely and uniformly on every compact subset of $\Omega$ 
also if $\alpha-2m>0$.

\item For any choice of $\alpha,m,l,t$ submitted to the convergence conditions above, the function ${}^t\mathcal{E}_{\alpha,m,l}(z,t)$ 
belongs to the space $\mathcal{M}^{l+1}_{\alpha,m}({}^t\rho_{t,l}^{-1},r(\alpha,m,l))$.

\item If $\alpha-l\not\equiv 2m\pmod{(q-1)}$, the matrix function $\mathcal{E}_{\alpha,m,l}(z,t)$ is identically zero.

\item If $\alpha-l\equiv 2m\pmod{(q-1)}$, $\alpha\geq (q+1)m$ and if $\mathcal{E}_{\alpha,m,l}$ converges, then $\mathcal{E}_{\alpha,m,l}$ is not identically zero in its domain of convergence.
\end{enumerate}
\end{Proposition}

\noindent\emph{Proof.} 
\noindent\emph{1.} We choose $\delta\in H\backslash \Gamma$ corresponding
to a couple $(c,d)\in A^2$ with $c,d$ relatively prime, and set $f_\delta=\mu_{\alpha,m}(\delta,z)^{-1}u^m(\delta(z))V_l (\delta)\in\mathcal{S}$.
We have
\begin{eqnarray*}
f_\delta(\gamma(z))&=&\mu_{\alpha,m}(\delta,\gamma(z))^{-1}u^m(\delta(\gamma(z)))V_l (\delta)\\
&=&\mu_{\alpha,m}(\gamma,z)\mu_{\alpha,m}(\delta\gamma,z)^{-1}u^m(\delta\gamma(z)))V_l (\delta),\\
&=&\mu_{\alpha,m}(\gamma,z)\mu_{\alpha,m}(\delta\gamma,z)^{-1}u^m(\delta\gamma(z)))V_l (\delta\gamma)\cdot\rho_{t,l}(\gamma)^{-1},\\
&=&\mu_{\alpha,m}(\gamma,z)\mu_{\alpha,m}(\delta',z)^{-1}u^m(\delta'(z))V_l (\delta')\cdot\rho_{t,l}(\gamma)^{-1},\\
&=&\mu_{\alpha,m}(\gamma,z)f_{\delta'}\cdot\rho_{t,l}(\gamma)^{-1},
\end{eqnarray*}
with $\delta'=\delta\gamma$ and $f_{\delta'}=\mu_{\alpha,m}(\delta',z)^{-1}u^m(\delta'(z))V_l (\delta')$, from which part 1 of the
proposition follows.

\medskip

\noindent\emph{2.}
Convergence and holomorphy are ensured by simple modifications of \cite[(5.5)]{Ge}, or by the arguments in \cite[Chapter 10]{GePu}.
More precisely, let us choose an integer $0\leq s\leq l$ and look at the component at the place $s$ $$\mathcal{E}_s(z,t)=\sum_{\delta\in H\backslash\Gamma}\mu_{\alpha,m}(\delta,z)^{-1}u(\delta(z))^m\ol{c^sd^{l-s}}$$ of the vector 
series $\mathcal{E}_{\alpha,m,l}$. Writing $\alpha=n(q-1)+2m+l'$ with $n$ non-negative integer and $l'\geq l$ we see, following \cite[pp. 304-305]{GePu} and taking into account 
$|u(\delta(z))|\leq|cz+d|^2/|z|_i$ ($|z|_i$ denotes, for $z\in \CC_\infty$, 
the infimum $\inf_{a\in K_\infty}\{|z-a|\}$), that
the term 
of the series $\mathcal{E}_s$:
$$\mu_{\alpha,m}(\delta,z)^{-1}u^m(\delta(z))\ol{c^sd^{l-s}}=(cz+d)^{-n(q-1)-l'-2m}u(\delta(z))^m\chi_t(c^sd^{l-s})$$
 (where $\delta$ corresponds to $(c,d)$)
has absolute value bounded from above by 
$$|z|_i^{-m}\left|\frac{\ol{c^sd^{l-s}}}{(cz+d)^{n(q-1)+l'}}\right|.$$
Taking into account the first part of the proposition, 
to check convergence, we can freely substitute $z$ with $z+a$ with $a\in A$ and we may assume, without loss of generality, that
$\deg_\theta z=\lambda\not\in\ZZ$. In this case, for all $c,d$, $|cz+d|=\max\{|cz|,|d|\}$.
Then, the series defining $\mathcal{E}_s$ can be decomposed as follows:
$$\mathcal{E}_s=\sum_{f_\delta\in H\backslash\Gamma}f_\delta=\left(\sideset{}{'}\sum_{|cz|<|d|}+\sideset{}{'}\sum_{|cz|>|d|}\right)\mu_{\alpha,m}(\delta,z)^{-1}u^m(\delta(z))\ol{c^sd^{l-s}}.$$
We now look for upper bounds for the absolute values of the terms of the series above separating the two cases in a way similar to that of Gerritzen and van der Put
in loc. cit. 

Assume first that $|cz|<|d|$, that is, $\deg_\theta c+\lambda< \deg_\theta d$. Then
$$\left|\frac{\ol{c^sd^{l-s}}}{(cz+d)^{n(q-1)+l'}}\right|\leq \kappa \max\{1,|t|\}^{l\deg_\theta d}|d|^{-n(q-1)-l'}\leq \kappa q^{\deg_\theta d(l\log^+_q|t|-n(q-1)-l')},$$
where $\kappa$ is a constant depending on $\lambda$, and the corresponding sub-series 
converges with the imposed conditions on the parameters, because $l\log^+_q|t|-n(q-1)-l'<0$.

If on the other side $|cz|>|d|$, that is, $\deg_\theta c+\lambda> \deg_\theta d,$ then
$$\left|\frac{\ol{c^sd^{l-s}}}{(cz+d)^{n(q-1)+l'}}\right|\leq\kappa' \max\{1,|t|\}^{l\deg_\theta d}|c|^{-n(q-1)-l'}\leq\kappa' q^{\deg_\theta c(l\log^+_q|t|-n(q-1)-l')},$$
with a constant $\kappa'$ depending on $\lambda$, again because $l\log^+_q|t|-n(q-1)-l'<0$. This completes the proof of the second part of the Proposition.

\medskip

\noindent\emph{3.} This property can be deduced from the proof of the second part because if $\log_q|t|<0$, then $|\chi_t(c^sd^{l-s})|\leq 1$. 

\medskip

\noindent\emph{4.} The property is obvious by the first part of the proposition, because 
$\mathcal{E}_{\alpha,m,l}=\sum_{f\in\mathcal{S}}f.$

\medskip

\noindent\emph{5.}
We consider $\gamma=\mathbf{Diag}(1,\lambda)$ with $\lambda\in\FF_q^\times$; the 
corresponding homography, multiplication by $\lambda^{-1}$, is equal to that defined by $\mathbf{Diag}(\lambda^{-1},1)$. Hence, we have:
\begin{eqnarray*}
\mathcal{E}_{\alpha,m,l}(\gamma(z))&=&\lambda^{\alpha-m}\mathcal{E}_{\alpha,m,l}(z)\cdot\mathbf{Diag}(1,\lambda^{-1},\ldots,\lambda^{-l})\\
&=&\lambda^{m}\mathcal{E}_{\alpha,m,l}(z)\cdot\mathbf{Diag}(\lambda^l,\lambda^{l-1},\ldots,1),
\end{eqnarray*}
from which it follows that $\mathcal{E}_{\alpha,m,l}$ is identically zero if $\alpha-l\not\equiv 2m\pmod{q-1}$.

\medskip

\noindent\emph{6.}
It is easy to modify the arguments in the proof of \cite[Proposition 10.5.2]{GePu}, where the case $l=0$ is handled.
Indeed, let us choose the value $t=0$ and consider any component of the vector $\mathcal{E}_{\alpha,m,l}|_{t=\theta}(\sqrt{\theta})$.

Just as in \cite{GePu}, the sum can be again decomposed into three terms $A,B,C$, submitted to the same estimates 
as on pp. 305-306 of loc. cit., from which we deduce right away that with the conditions above on $\alpha,m$, the function $\mathcal{E}_{\alpha,m,l}$
is not identically zero.\CVD

Let $\alpha,m,l$ be non-negative integers such that $\alpha-2m>l$ and $\alpha-l\equiv 2m\pmod{(q-1)}$. We have functions:
\begin{eqnarray*}
\mathcal{E}_{\alpha,m,l}:\Omega&\rightarrow&\mathbf{Mat}_{1\times l+1}(\mathcal{R}_{<r}),\\
\Phi_l:\Omega&\rightarrow&\mathbf{Mat}_{l+1\times 1}(\mathcal{R}_{\infty}),
\end{eqnarray*}
with $r=r(\alpha,m,l)$ as in Proposition \ref{mainproppoincare}, and 
${}^t\mathcal{E}_{\alpha,m,l}\in\mathcal{M}^{l+1}_{\alpha,m}({}^t\rho_{t,l}^{-1},r)$, $\Phi_l\in\mathcal{M}^{l+1}_{-l,0}(\rho_{t,l},\infty)$. Therefore,
after Proposition \ref{equationf}, the functions $$\mathcal{G}_{\alpha,m,l,k}=(\tau^k\mathcal{E}_{\alpha,m,l})\cdot\Phi_l=\mathcal{E}_{q^k\alpha,m,l}\cdot\Phi_l:\Omega\rightarrow\TT_{<r}$$ satisfy $\mathcal{G}_{\alpha,m,l,k}\in M^!_{q^k\alpha-l,m}\otimes\TT_r$.

\subsubsection{A special case: vectorial Eisenstein series}

After Proposition \ref{mainproppoincare},
if $\alpha>0,l\geq 0$ and $\alpha\equiv l\pmod{q-1}$, then $\mathcal{E}_{\alpha,0,l}\neq0$. We call these series {\em deformations of vectorial Eisenstein series}), and we focus especially on the case $l=1$. We introduce the following formal series ``value of $L$-series":
$$L(\chi_t^l,\alpha)=\sum_{a\in A^+}a^{-\alpha}\chi_t(a)^l\in K_\infty[[t]],$$ where $A^+$ denotes the set of monic polynomials of $A$.
After the inequality $$|\ol{a}^la^{-\alpha}|\leq q^{\deg_\theta a(l\log_q^+|t|-\alpha)},$$
the series $L(\chi_t^l,\alpha)$ converges for all $t$ such that $\log_q^+|t|<\alpha/l$ if $l\neq0$ (otherwise, there is no dependence on $t$). In particular, if $\alpha>l+1$, the series converges at $t=\theta$
to the Carlitz-Goss zeta value $\zeta(\alpha-l)=\sum_{a\in A^+}a^{\alpha-l}$ and if $\alpha>0$, the series converges at $t\in\FF_q^{\text{alg.}}$ 
to the value at the integer $\alpha$ of the $L$-series associated to a Dirichlet character.
Moreover, we have the following obvious relation, which helps us to extend the domain of definition of $L(\chi_t^l,\alpha)$ (analog of ``integration by parts"):
\begin{equation}\label{tauzeta}
\tau L(\chi_t^l,\alpha)=L(\chi_t^l,q\alpha).
\end{equation}

\begin{Lemme}\label{interpret} With $\alpha,l$ such that $\alpha\equiv l\pmod{q-1}$ and $\alpha\geq l$, the following identity holds:
$$\mathcal{E}_{\alpha,0,l}(z,t)=L(\chi_t^l,\alpha)^{-1}\sideset{}{'}\sum_{c,d}(cz+d)^{-\alpha} V_l(c,d),$$
and $\mathcal{E}_{\alpha,0,l}$ is not identically zero.
\end{Lemme}
\noindent\emph{Proof.}
We recall the notation $$V_l(c,d)=(\ol{c}^l,\ol{c}^{l-1}\ol{d},\ldots,\ol{d}^l)\in\mathbf{Mat}_{1\times l+1}(\FF_q[t]).$$ We have
\begin{eqnarray*}
\sideset{}{'}\sum_{c,d}(cz+d)^{-\alpha} V_l(c,d)&=& \sum_{(c',d')=1}\sum_{a\in A^+}a^{-\alpha}(c'z+d')^{-\alpha} V_l(ac',ad')\\
&=&L(\chi_t^l,\alpha)\mathcal{E}_{\alpha,0,l}(z,t),
\end{eqnarray*}
where the first sum is over couples of $A^2$ distinct from $(0,0)$, while the second sum is over the couples $(c',d')$ of relatively prime
elements of $A^2$. Non vanishing of the function follows from Proposition \ref{mainproppoincare}.\CVD

In particular, if $l=0$, we obtain classical Eisenstein series up to a factor of proportionality:
$$\mathcal{E}_{\alpha,0,0}(z,t)=L(1,\alpha)^{-1}\sideset{}{'}\sum_{c,d}(cz+d)^{-\alpha}=\zeta(\alpha)^{-1}\sideset{}{'}\sum_{c,d}(cz+d)^{-\alpha}.$$

\section{Proof of the main results}

Following Gekeler \cite[Section 3]{Ge}, we recall that for all $\alpha>0$ there exists a polynomial $G_\alpha(u)\in \CC_\infty[u]$, called the 
{\em $\alpha$-th Goss polynomial},
such that, for all $z\in\Omega$, $G_{\alpha}(u(z))$ equals the sum of the convergent series
$$\widetilde{\pi}^{-\alpha}\sum_{a\in A}\frac{1}{(z+a)^\alpha}.$$

Several properties of these polynomials are collected in \cite[Proposition (3.4)]{Ge}. Here, we will need that
for all $\alpha$, $G_\alpha$ is of type $\alpha$ as a formal series of $\CC_\infty[[u]]$. Namely:
$$G_\alpha(\lambda u)=\lambda^\alpha G_\alpha(u),\quad \text{ for all }\lambda\in\FF_q.$$

We also recall, for $a\in A$, the function
$$u_a(z):=u(az)=e_\carlitz(\widetilde{\pi}az)^{-1}=u^{|a|}f_a(u),$$
where $f_a\in A[[u]]$ is the {\em $a$-th inverse cyclotomic polynomial} defined in \cite[(4.6)]{Ge}. Obviously,
we have
$$u_{\lambda a}=\lambda^{-1}u_a\quad \text{ for all }\lambda\in\FF_q^\times.$$

To continue this section, we will state and prove three auxiliary lemmas.

\begin{Lemme}\label{primolemma}
Let $\alpha$ be a positive integer such that $\alpha\equiv 1\pmod{q-1}$.
We have, for all $t\in \CC_\infty$ such that $|t|<1$ and $z\in\Omega$, convergence of both the series below, and equality:
$$\sideset{}{'}\sum_{c,d\in A}\frac{\ol{c}}{(cz+d)^{\alpha}}=-\widetilde{\pi}^\alpha\sum_{c\in A^+}\ol{c}G_{\alpha}(u_c(z)).$$ 
\end{Lemme}
\noindent\emph{Proof.} Convergence features are easy to deduce from Proposition \ref{mainproppoincare}. 
Indeed, for $l=1$ we have convergence if $\log_q^+|t|<r(\alpha,m,l)=\alpha-1$, that is, $\max\{1,|t|\}\leq q^{\alpha-1}$ if $\alpha>1$
and we have convergence, for $\alpha=1$, for $|t|<1$. In all cases, convergence holds for $|t|<1$. 

We then compute:
\begin{eqnarray*}
\sideset{}{'}\sum_{c,d}\frac{\ol{c}}{(cz+d)^\alpha}&=&\sum_{c\neq0}\ol{c}\sum_{d\in A}\frac{1}{(cz+d)^\alpha}\\
&=&\widetilde{\pi}^\alpha\sum_{c\neq 0}\ol{c}\sum_{d\in A}\frac{1}{(c\widetilde{\pi}z+d\widetilde{\pi})^\alpha}\\
&=&\widetilde{\pi}^\alpha\sum_{c\neq 0}\ol{c}G_\alpha(u_c)\\
&=&\widetilde{\pi}^\alpha\sum_{c\in A^+}\ol{c}G_\alpha(u_c)\sum_{\lambda\in\FF_q^\times}\lambda^{1-\alpha}\\
&=&-\widetilde{\pi}^\alpha\sum_{c\in A^+}\ol{c}G_\alpha(u_c).
\end{eqnarray*}
\CVD

\begin{Lemme}\label{lemmelimit1} Let $\alpha>0$ be an integer such that $\alpha\equiv1\pmod{q-1}$. For all $t\in \CC_\infty$ such that $|t|<1$, we have
$$\lim_{|z|_i=|z|\to\infty}\bsb{d}_1(z)\sideset{}{'}\sum_{c,d}\frac{\ol{c}}{(cz+d)^\alpha}=0.$$
\end{Lemme}
\noindent\emph{Proof.} 
We recall from \cite{archiv} the series expansion
$$\bsb{d}_1(z)=\frac{\widetilde{\pi}}{s_{\carlitz}(t)}\bsb{s}_2(z)=\frac{\widetilde{\pi}}{s_{\carlitz}(t)}\sum_{n\geq0}e_{\Lambda_z}\left(\frac{z}{\theta^{n+1}}\right)t^n,$$
converging for all $t$ such that $|t|<q$ and all $z\in\Omega$.

By a simple modification of the proof of \cite[Lemma 5.9 p. 286]{gekeler:compositio}, we have
$$\lim_{|z|_i=|z|\to\infty}u(z)t^ne_{\Lambda_z}(z/\theta^{n+1})^q=0$$ uniformly in $n>0$, for all $t$ such that
$|t|\leq q$.

Moreover, it is easy to show that
$$\lim_{|z|_i=|z|\to\infty} u(z)e_{\Lambda_z}(z/\theta)^q=\widetilde{\pi}^{-q}\lim_{|z|_i=|z|\to\infty}e^{q}_{\carlitz}(\widetilde{\pi}z/\theta)/e_{\carlitz}(\widetilde{\pi}z)=1.$$
This suffices to show that
$$\lim_{|z|_i=|z|\to\infty}\bsb{d}_1(z)G_{\alpha}(u_c(z))=0$$
uniformly for $c\in A^+$, for all $t$ such that $|t|<q$. The lemma then follows from the application of Lemma \ref{primolemma}.\CVD 

\begin{Lemme}\label{lemmelimit2} Let $\alpha>0$ be an integer such that $\alpha\equiv1\pmod{q-1}$. For all $t\in \CC_\infty$ such that $|t|<1$, we have
$$\lim_{|z|_i=|z|\to\infty}\sideset{}{'}\sum_{c,d}\frac{\ol{d}}{(cz+d)^\alpha}=-L(\chi_t,\alpha).$$
\end{Lemme}
\noindent\emph{Proof.} It suffices to 
show that
$$\lim_{|z|_i=|z|\to\infty}\sum_{c\neq0}\sum_{d\in A}\frac{\ol{d}}{(cz+d)^{\alpha}}=0.$$
We assume, as we can, that $z\in\Omega$ is chosen so that, for all $(c,d)\in A\setminus\{(0,0)\}$,
$|cz|\neq|d|$. We then have, for $c\neq 0$:
\begin{eqnarray*}
\sum_{d\in A}\frac{\ol{d}}{(cz+d)^{\alpha}}&=&\left(\sum_{|d|>|cz|}+\sum_{|d|<|cz|}\right)\frac{\ol{d}}{(cz+d)^{\alpha}}.
\end{eqnarray*}
Now, if $|d|>|cz|,$ we have, for $|t|<1$:
$$
\left|\frac{\chi_t(d)}{(cz+d)^{\alpha}}\right|\leq\left|\frac{\chi_t(d)}{d^{\alpha}}\right|\leq|d|^{-\alpha}\leq |cz|^{-\alpha}.
$$
If $|d|<|cz|,$ we have, again for $|t|<1$:
$$
\left|\frac{\chi_t(d)}{(cz+d)^{\alpha}}\right|\leq\left|\frac{\chi_t(d)}{(cz)^{\alpha}}\right|\leq|cz|^{-\alpha}.$$
Therefore, for $c\neq0$,
$$\left|\sum_{d\in A}\frac{\ol{d}}{(cz+d)^{\alpha}}\right|\leq|cz|^{-\alpha}.$$
This implies that 
$$\left|\sum_{c\neq0}\sum_{d\in A}\frac{\ol{d}}{(cz+d)^{\alpha}}\right|\leq|z|^{-\alpha},$$
from which the Lemma follows.\CVD

The next step is to prove the following Proposition.

\begin{Proposition}\label{interpretationvk}
Let $\alpha$ be positive, such that $\alpha\equiv 1\pmod{q-1}$. Then,
the sequence
$$(\mathcal{G}_{\alpha,0,1,k})_{k\in\ZZ}$$
is generic for the difference field $(\mathcal{K},\tau)$.
Moreover, if $\alpha\leq q(q-1)$,
then:
$$\mathcal{G}_{\alpha,0,1,0}=-E_{\alpha-1},$$
where $E_{\alpha-1}$ is the normalised Eisenstein series of weight $\alpha-1$.
\end{Proposition}

\noindent\emph{Proof.}
The components of $\Phi_1$ are $\FF_q(t)$-linearly independent (Lemma \ref{componentslinearlyindep}). The $\FF_q(t)$-linear 
independence of the components of $\mathcal{E}_{\alpha,0,1}$ follows from analysing the behaviour at $u=0$ described by Lemmas
\ref{lemmelimit1} and \ref{lemmelimit2}. This means that the sequence $(\mathcal{G}_{\alpha,0,1,k})_{k\in\ZZ}$ is generic hence proving the
first part of the proposition.

According to Lemma \ref{interpret}, we need, to finish the proof of the proposition, to compute the sum of the series:
$$F_\alpha(z):=\bsb{d}_1(z)\sideset{}{'}\sum_{c,d}\frac{\ol{c}}{(cz+d)^\alpha}+\bsb{d}_2(z)\sideset{}{'}\sum_{c,d}\frac{\ol{d}}{(cz+d)^\alpha},$$
which converges in $\Omega$, noticing that this yields the case $\alpha<q(q-1)$ in Theorem \ref{scalarprod}. 

After (\ref{uexpd2}), we have that for all $t$ with $|t|<1$,
$\lim_{|z|_i=|z|\to\infty}\bsb{d}_2(z)=1$. From Lemmas \ref{lemmelimit1} and \ref{lemmelimit2},
$$\lim_{|z|_i=|z|\to\infty}F_\alpha(z)=-L(\chi_t,\alpha).$$

In particular, $F_\alpha(z)$ is a modular form of $M_{\alpha-1,0}\otimes\TT_{<q}$. 
Since for the selected values of $\alpha$, $M_{\alpha-1,0}=\langle E_{\alpha-1}\rangle$, 
we obtain that $F_\alpha=-L(\chi_t,\alpha)E_{\alpha-1}$. After Lemma \ref{interpret},
the proposition follows.\CVD

\subsection{Proofs of the main theorems} 

We prove Theorem \ref{firsttheorem} here and Theorems \ref{charroots}, \ref{scalarprod} are simple consequences of it. We will also 
deduce Corollaries  \ref{corollairezeta11}, \ref{firstcorollary}. At the end of the subsection, there is a proof of Theorem \ref{propositionxk}.

\medskip

\noindent\emph{Proofs of Theorem \ref{firsttheorem} and Corollaries \ref{corollairezeta11}, \ref{firstcorollary}.}
Let us write:
$$\mathcal{E}=\mathcal{E}_{1,0,1}=L(\chi_t,1)^{-1}\sideset{}{'}\sum_{c,d\in A}\left(\frac{\ol{c}}{cz+d},
\frac{\ol{d}}{cz+d}\right),\quad \mathcal{F}=\Phi_1=\binom{\bsb{d}_1}{\bsb{d}_2},\quad \mathcal{G}_k= \mathcal{G}_{1,0,1,k}=(\tau^k\mathcal{E})\cdot\mathcal{F}.$$
With the notations of the introduction, we have
$$\mathcal{E}=L(\chi_t,1)^{-1}(\bsb{e}_1,\bsb{e}_2).$$

We know by (\ref{equexpansion}) that the entries of $\mathcal{F}$ span the $\FF_q(t)$-vector space of solutions in $\mathcal{K}$ of
the $\tau$-difference equation
$Lx=0,$
where $L=L_1$ is the operator defined in (\ref{ell1}),
and Proposition \ref{equationf} implies that $L(\mathcal{G}_k)\equiv 0$. Combining Lemma \ref{r=s} with Proposition \ref{prarchiv}
we finally obtain Theorem \ref{firsttheorem}.\CVD

Comparing the coefficients of $u$ in the $u$-expansions of both sides of the identity
$$L(\chi_t,1)^{-1}\bsb{e}_1=-h\tau (s_\carlitz\bsb{d}_2),$$ obtained from Theorem \ref{firsttheorem} (with the help of Proposition \ref{prarchiv}), and using Lemma \ref{primolemma}
one deduces Corollary \ref{corollairezeta11}. Replacing $L(\chi_t,1)$ by $-\widetilde{\pi}/(\tau s_\carlitz)$ in the latter identity and using again
Lemma \ref{primolemma} and the definition of $\bsb{E}$ then yields Corollary \ref{firstcorollary}.\CVD

\noindent\emph{Proof of Theorems \ref{charroots} and \ref{scalarprod}.} 
For $k=0,1$, Theorems \ref{charroots} and \ref{scalarprod} agree with Proposition \ref{interpretationvk}. Since $\mathcal{G}=(\mathcal{G}_k)_{k\in\ZZ}$
satisfies $$L_1(\mathcal{G})=0$$ where $L_1$ is the operator defined in (\ref{ell1}) and by definition,
$L_1((g^\star_k)_{k\in\ZZ})=0$, the two sequences $\mathcal{G}$ and $(g^\star_k)_{k\in\ZZ}$ have the same 
initial data, so they are equal. The $\tau$-linearised recurrent relations (\ref{secondrecurrence}) are easy to 
obtain by computing explicitly the operator $L'=A_0\tau^0+A_1'\tau^{-1}+A_2'\tau^{-2}\in \mathcal{K}[\tau^{-1}]$ such that the entries of $\mathcal{E}$ span the $\FF_q(t)$-vector space 
of solutions in $\mathcal{K}$ of $L'x=0.$ We find $$L'=-\tau^0+g^{1/q}\tau^{-1}+\Delta^{1/q^2}(1-\theta^{1/q})\tau^{-2}.$$
\CVD

\noindent\emph{Proof of Theorem \ref{propositionxk}.}
We recall that after \cite{archiv,BP3}, for all $k\geq0$, $\bsb{E}$ is a {\em deformation of Drinfeld quasi-modular form} 
of weight $(q^k,1)$ and type $1$
and the function $\phi_k(z):=\boldsymbol{E}^{(k)}(z,\theta)$ is a well defined
Drinfeld quasi-modular form in the space $\widetilde{M}_{q^k+1,1}^{\leq 1}$.
By Corollary \ref{firstcorollary}, $$\bsb{E}^{(k)}=u^{q^k}+\cdots,$$
and again by \cite[Theorem 1.2]{BP2} $\bsb{E}^{(k)}(z,\theta)$ is normalised, extremal, therefore 
proportional to $x_k$ for all $k$. By \cite[Proposition 2.3]{BP2},
$$x_k=(-1)^{k+1}L_ku^{q^k}+\cdots,$$ where $L_k=[k][k-1]\cdots[1]$ if $k>0$ and
$L_0=1$. From Corollary \ref{firstcorollary}, for all $k\geq 0$, $\bsb{E}^{(k)}\in \FF_q[t,\theta][[u]]$.
Therefore, $E_k\in A[[u]]$.\CVD

\noindent\emph{Remarks. 1.} As we already mentioned, from (\ref{eqP101}) we find $\mathcal{E}_{1,0,1}\cdot\Phi_1=\mathcal{G}_0=-1$, which is
our deformation of Legendre's identity (\ref{dethPsi}).

\medskip

\noindent\emph{2.} Corollary \ref{corollairezeta11} implies that for all $k\geq 0$, $\tau^kL(\chi_t,1)=-\widetilde{\pi}^{q^k}/(\tau^ks_{\carlitz}^{(1)})$.
Since $\tau^ks_\carlitz^{(1)}=(t-\theta^{q^k})\cdots(t-\theta^q)s_\carlitz^{(1)}$, we obtain, for $k>0$, the well known formulas for Carlitz-Goss' zeta values
$$\zeta(q^k-1)=(-1)^k\frac{\widetilde{\pi}^{q^k-1}}{[k][k-1]\cdots[1]}.$$ 

\medskip

\noindent\emph{3.} 
Proposition \ref{equationf} states the existence of 
$L'=A_0\tau^0+A_1'\tau^{-1}+A_2'\tau^{-2}\in \mathcal{K}[\tau^{-1}]$ such that the entries of $\mathcal{E}$ span the $\FF_q(t)$-vector space 
of solutions in $\mathcal{K}$ of $L'x=0.$
From Theorem  \ref{firsttheorem} we said it is easy to deduce that $$L'=-\tau^0+g^{1/q}\tau^{-1}+\Delta^{1/q^2}(1-\theta^{1/q})\tau^{-2}.$$
This is the $\tau^{-1}$-form of the {\em adjoint} of $L$ of \cite[Goss, Section 1.7]{Go2}, denoted by $L^*$ there. Keeping the 
notations of Goss, we then have the $\tau$-form of the adjoint, $L^{\text{ad}}=\tau^2L'\in\mathcal{K}[\tau]$:
$$L^{\text{ad}}=-\tau^2+g^q\tau+\Delta(t-\theta^q)\tau^0.$$
The fact that, simultaneously, $L(\Phi_1)=0$ and $L^{\text{ad}}(\mathcal{E}_{1,0,1})=0$ ($\alpha=1$) is a peculiar phenomenon which
does not seem to hold for general values of $\alpha$. It would be interesting to understand when this occurs.

\medskip

\section{Computing $u$-expansions\label{computing}}

Let $\mu$ be an element of $\CC_\infty$ and let us consider the function:

$$s_{\carlitz,\mu}(t):=\sum_{i=0}^\infty e_\carlitz\left(\frac{\mu}{\theta^{i+1}}\right)t^i.$$
The function $\mu\mapsto s_{\carlitz,\mu}$ is well defined with image in $\TT_{<q}$.

By \cite[Equation (10) p. 220]{Bourbaki} we have the functional equation:

\begin{equation}\label{eqsu}
s_{\carlitz,\mu}^{(1)}(t)=(t-\theta)s_{\carlitz,\mu}(t)+e_{\carlitz}(\mu).
\end{equation}
For fixed $\mu,$ the function $s_{\carlitz,\mu}(t)$ has a simple pole in $t=\theta$ with residue $-\mu$.
We point out that $s_\carlitz=s_{\carlitz,\widetilde{\pi}}$. 

We now consider the function $$F^\star:\CC_\infty\rightarrow\TT_{<q}$$ defined by $F^\star(z)=s_{\carlitz,\widetilde{\pi}z}(t)\in\TT_{<q}$ 
(so that $F^\star(1)=s_\carlitz$ and $F^\star\in\mathcal{R}_{<q}$) and
the function $F:\Omega\rightarrow \TT_{<q}$ defined by $F(z)=F^\star(z)/s_\carlitz$. We have $F\in\mathcal{R}_\infty$
and we can write:
\begin{equation}\label{valueFtheta}
F(z)|_{t\mapsto\theta}=z.
\end{equation}
We have the functional
equations
\begin{equation}\label{eqfunctsu}
F^{(1)}=F+\frac{1}{(t-\theta)us_\carlitz(t)},\quad F^\star{}^{(1)}=(t-\theta)F^\star+\frac{1}{u}.
\end{equation}

In the next two propositions, we introduce the functions $\psi,\bsb{d}_3,\psi^\star,\bsb{d}_3^\star$.
In fact, we set $\psi^\star=s_\carlitz^{(1)}\psi$ and $\bsb{d}_3^\star=s_\carlitz^{(1)} \bsb{d}_3 $ so
that we only need to define $\psi$ and $\bsb{d}_3$, but we will discuss properties of
all the four functions.

\begin{Proposition}\label{propositionfunctionpsi}
Let us define the function $\psi=R\bsb{d}_2+R^{(1)}(\bsb{d}_2-g\bsb{d}_2^{(1)})$ with $R=1/((t-\theta)us_\carlitz)=1/(us_\carlitz^{(1)})$ and let 
$\psi^\star$ be the function $$s_\carlitz^{(1)}\psi=\frac{\bsb{d}_2}{u}+\frac{\bsb{d}_2-g\bsb{d}_2^{(1)}}{(t-\theta^q)u^q}=\frac{\bsb{d}_2}{u}+\frac{\bsb{d}_2^{(2)}\Delta}{u^q}.$$ We have the following properties.
\begin{enumerate}
\item The function $\psi$ belongs to $\mathcal{R}_{<q^q}$.
\item The function $\psi^\star$ can be identified, for $|u|,|t|$ small, with the sum 
of a converging $u$-expansion
$$\psi^\star\in u^{q-2}\FF_q[t,\theta][[u^{q-1}]]$$ of type $-1$.
\item The first few terms of the $u$-expansion of $\psi^\star$ read as follows:
\begin{equation}\label{expansionpsi}
u^{q-2}(\theta-t+u^{(q-1)(q-2)}+(\theta-\theta^q)u^{(q-1)^2}+\cdots).
\end{equation}
\item We have, for all $u$ with $|u|$ small enough, 
$$\lim_{t\to\theta}\psi^\star=\frac{1}{u}+\frac{Eg+h}{(\theta-\theta^q)hu^q}.$$
\item If $q\neq 2$, we have $\lim_{u\to 0}\psi^\star=0$, while if $q=2$, we have $\lim_{u\to 0}\psi^\star=1+\theta-t$.
\end{enumerate}
\end{Proposition}

\noindent\emph{Proof. 1.} This is clear as $\bsb{d}_2,\bsb{d}_2^{(2)}$ belong to $\mathcal{R}_{<q^q}$
as well as $R=1/(us_\carlitz^{(1)})$.

\medskip

\noindent\emph{2.} Writing $v=u^{q-1}$, we have by (\ref{d2inu}):
$$\bsb{d}_2=1+(\theta-t)v+(\theta-t)v^{q^2-q+1}+\cdots$$
and we have the series expansion
$$g=1+(\theta-\theta^q)v+(\theta-\theta^q)v^{q^2-q+1}+\cdots,$$
that can be obtained with \cite[Corollary (10.11) and formula for $U_1$ on p. 691]{Ge}.

Substituting into the definition of $\psi$, we see, from $\bsb{d}_2\in\FF_q[t,\theta][[v]]$ that $$\psi^\star\in u^{q-2}\FF_q[t,\theta][[v]].$$
Moreover, we know that $\Delta,\bsb{d}_2,\bsb{d}_2^{(2)}$ are of type $0$, and it is obvious that
$R^{(k)}$ is of type $-1$ for all $k$.

\medskip

\noindent\emph{3.} 
Explicitly, we compute step by step:
\begin{eqnarray*}
\bsb{d}_2^{(1)}&=&1+(\theta^q-t)v^q+(\theta^q-t)v^{q(q^2-q+1)}+\cdots\\
g\bsb{d}_2^{(1)}&=&1+(\theta-\theta^q)v+(\theta^q-t)v^q+(\theta^q-\theta)(t-\theta^q)v^{q+1}+(\theta-\theta^q) v^{q^2-q+1}+\cdots\\
\bsb{d}_2-g\bsb{d}_2^{(1)}&=&(\theta^q - t)(v - v^q + (\theta^q - \theta)v^{q + 1} + v^{q^2 - q + 1}+\cdots)\\
\bsb{d}_2+\frac{\bsb{d}_2-g\bsb{d}_2^{(1)}}{(t-\theta^q)v}&=&(\theta-t)v+v^{q-1}+(\theta-\theta^q)v^q-v^{q^2-q}+(\theta-t)v^{q^2-q+1}+\cdots\\
\psi^\star&=&\frac{1}{u}\left(\bsb{d}_2+\frac{\bsb{d}_2-g\bsb{d}_2^{(1)}}{(t-\theta^q)v}\right)\\
&=&u^{-1}\{(\theta-t)v+v^{q-1}+(\theta-\theta^q)v^q-v^{q^2-q}+(\theta-t)v^{q^2-q+1}+\cdots\},
\end{eqnarray*}
which gives (\ref{expansionpsi}) and all the properties of $\psi$ claimed by the statement of the proposition.

\medskip

\noindent\emph{4.} It suffices to use the definition of $\psi^\star$, $\bsb{d}_2^{(2)}=\frac{\bsb{d}_2-g\bsb{d}_2^{(1)}}{(t-\theta^q)\Delta}$, and the identities $\bsb{E}=-h\bsb{d}_2^{(1)}$
and $\bsb{E}(\theta)=E,\bsb{d}_2(\theta)=1$.

\medskip

\noindent\emph{5.} This follows directly from (\ref{expansionpsi}).\CVD

\begin{Proposition} \label{lemmed1}
We define the functions $\bsb{d}_3=\bsb{d}_1-\bsb{d}_2F$ and $\bsb{d}_3^\star=s_\carlitz^{(1)}\bsb{d}_3$. The following properties hold.

\begin{enumerate}
\item We have that $\bsb{d}_3\in\mathcal{R}_{\infty}$.
\item The function $\bsb{d}_3$ is solution of the 
non-homoge\-neous $\tau$-difference equation:
\begin{equation}\label{eqdelta2}
\bsb{d}_3=(t-\theta^q)\Delta\bsb{d}_3^{(2)}+g\bsb{d}_3^{(1)}+\psi,
\end{equation}
\item The function $\bsb{d}_3^\star$ can be identified, for $|u|,|t|$ small, with the sum 
of a converging series:
$$\bsb{d}_3^\star\in\FF_q[t,\theta][[u]]$$ of type $-1$.
The $u$-expansion of $\bsb{d}_3^\star$ begins, for $q\neq2$, with the following terms:
\begin{equation}\label{qdifferent2}
-u^{q-2}(t-\theta+(t-\theta)u^{q(q-1)^2}+\cdots).
\end{equation}
If $q=2$, the $u$-expansion of $\bsb{d}_3^\star$ begins with the following terms:
\begin{equation}\label{qequal2}t+\theta+(1+t+\theta)u^2+\cdots.
\end{equation}
\item We have the limit $\lim_{t\to\theta}\bsb{d}_3=0$ for all $z\in\Omega$ and $\bsb{d}_3$ is the only 
solution of (\ref{eqdelta2}) with this property.
\end{enumerate}
\end{Proposition}
\noindent\emph{Proof. 1.} We have seen that $F,\bsb{d}_2$ are in $\mathcal{R}_\infty$,
so that the property follows for $\bsb{d}_3$.

\medskip\noindent\emph{2.} According to (\ref{eqsu}), we get:
\begin{eqnarray*}
\bsb{d}_1^{(1)}&=&\bsb{d}_2^{(1)}F^{(1)}+\bsb{d}_3^{(1)}\\
&=&\bsb{d}_2^{(1)}(F+R)+\bsb{d}_3^{(1)},\\
\bsb{d}_1^{(2)}&=&\bsb{d}_2^{(2)}F^{(2)}+\bsb{d}_3^{(2)}\\
&=&\bsb{d}_2^{(2)}(F+R+R^{(1)})+\bsb{d}_3^{(2)}.
\end{eqnarray*}
Let $L$ be the operator $L_1$ defined in (\ref{ell1}).
By (\ref{equexpansion}), we have 
$L\bsb{d}_1=0$. Explicitly,
$$F\bsb{d}_2+\bsb{d}_3=(t-\theta^q)\Delta((F+R+R^{(1)})\bsb{d}_2^{(2)}+\bsb{d}_3^{(2)})+g((F+R)\bsb{d}_2^{(1)}+\bsb{d}_3^{(1)}).$$
But again by (\ref{equexpansion}), $L\bsb{d}_2=0$ and we see that all the coefficients of $F$ in the identity above give contribution $0$
(alternatively, we can apply Lemma \ref{lemmatranscendence} and the fact that $\bsb{d}_2$ is a formal power series in $u$).
In other words,
$$L\bsb{d}_3+(t-\theta^q)\Delta(R+R^{(1)})\bsb{d}_2^{(2)}+Rg\bsb{d}_2^{(1)}=0.$$
Eliminating $\bsb{d}_2^{(2)}$ with (\ref{equexpansion})
in the above expression yields
\begin{equation}\label{delta2}
L\bsb{d}_3+\psi=0,\end{equation}
that is, (\ref{eqdelta2}).

\medskip\noindent\emph{3, 4, 5.} 
We proceed as in \cite{archiv}, where we computed the $u$-expansion (\ref{d2inu}).
We first look at the case $q\neq2$ and then, we consider the case $q=2$, more involved.
We begin by showing that for $q\neq2$ equation (\ref{eqdelta2}) has an unique solution
${Y}$ which can be expanded as a formal series in powers of $u$, with the property 
that ${Y}|_{t=\theta}=0$. Then, we show that $\bsb{d}_3={Y}$.

Let $f$ be a formal series in non-negative powers of $u$ with coefficients, say, 
in $\TT_{<r}$,
$$f=\sum_ic_iu^i,$$ and let us consider the truncation $$[f]_n=\sum_{i\leq q^{n}-1}c_iu^i$$ of the series $f$ to the order $q^n$, with $n\geq 0$
(do not mix up with the truncation in powers of $t$ also used in this paper). By convention, we also set $[f]_n=c_0$ for $n<0$.
We have, for series $f,g$, the following simple identities:
\begin{enumerate}
\item $[f+g]_n=[f]_n+[g]_n$,
\item $[fg]_n=[[f]_n[g]_n]_n$,
\item $[f^{(1)}]_n=[f]_{n-1}^{(1)}$.
\end{enumerate}
For all $n\geq 2$ and any series ${Y}=\sum_{i\geq 0}c_iu^i$ solution of (\ref{eqdelta2}),
$$[{Y}]_n=(t-\theta^q)[[\Delta]_n[{Y}]_{n-2}^{(2)}]_n+[[g]_n[{Y}]_{n-1}^{(1)}]_n+[\psi]_n.$$

Hence, if ${Y}$ exists, the whole collection of its coefficients is uniquely determined by $[{Y}]_1$, and the integrality of the coefficients 
of $[Y]_n$ follows from the same property for $[Y]_1$.
We recall now that we are assuming that $q\neq 2$. In this case, $\psi$ vanishes at $u=0$ (Proposition \ref{propositionfunctionpsi}) and for $n=1$, we find:
$$[{Y}]_1=[\psi]_1.$$
This means that there exists one and only one solution of (\ref{eqdelta2}) for $q\neq 2$ which is a series of powers of $u$, with the additional
property that it vanishes at $u=0$.

Now, we need to show that ${Y}$ is the function we are looking for, but this is a simple task.
The set of solutions in $\mathcal{R}_\infty$ of (\ref{eqdelta2}) is the 
translated of $\FF_q(t)$-vector space:
$$\FF_q(t)\bsb{d}_1+\FF_q(t)\bsb{d}_2+{Y}.$$
Since $\bsb{d}_1=\bsb{d}_2F+\bsb{d}_3$ and
$\bsb{d}_1|_{t=\theta}=F|_{t=\theta}=z$ and $\bsb{d}_2|_{t=\theta}=1$, we have $\bsb{d}_3|_{t=\theta}=0$ and we see that 
$$\bsb{d}_3={Y}.$$
The $u$-expansion (\ref{qdifferent2}) can be checked
after explicit computation.

Also, by induction, we may verify that $\bsb{d}_3,\bsb{d}_3^\star$ have type $-1$ and that
$$\bsb{d}_3^\star\in u^{q-2}\FF_q[t,\theta][[v]].$$

Let us now consider the case $q=2$, in which types are trivial and $u=v$. Here, $\psi$ does not vanish at $u=0$ and this case has to be handled in slightly 
different way. In this case, we have, returning to the unknown series ${Y}$, the identities:
$$[{Y}]_1=(t-\theta^q)[[\Delta]_1[{Y}]_{0}^{(2)}]_1+[[g]_1[{Y}]_{0}^{(1)}]_1+[\psi]_1$$
and 
$$[{Y}]_0=[[g]_0[{Y}]_{0}^{(1)}]_0+[\psi]_0.$$
Now, the truncations $[\Delta]_1$ and $[g]_1$ are easy to compute:
\begin{eqnarray*}
{[}\Delta{]}_1&=&-u,\\
{[}g{]}_1&=&1+(\theta^2+\theta)u
\end{eqnarray*}
By (\ref{expansionpsi}), $[\psi]_1$ is:
$$s_\carlitz^{-1}\left(1+\frac{1}{t+\theta}+(\theta+\theta^2)u\right).$$
Hence, the constant term $c_0=[{Y}]_0$ satisfies a $\tau$-difference ``Artin-Schreier" equation:
$$c_0=c_0^{(1)}+s_\carlitz^{-1}\left(1+\frac{1}{t+\theta}\right)$$
whose set of solutions is $\{s_\carlitz^{-1}+\lambda\}$ with $\lambda\in\FF_q(t)$ and we are reduced to calculate $\lambda$
corresponding to our function $\bsb{d}_3$, which satisfies $\bsb{d}_3|_{t=\theta}=0$. We deduce that 
$\lim_{t\to\theta}c_0=0.$ Therefore, $\lambda=0$ and after some computations, we find (\ref{qequal2}).
The reader can verify that all the properties of the proposition have been checked.\CVD

\noindent\emph{Remark.} It can be proved
that $F$ is, up to multiplication by a factor in $\FF_q(t)$, the only function for which we can write $\bsb{d}_1=\bsb{d}_2F+\bsb{d}_3$, with $\bsb{d}_2,\bsb{d}_3$
formal power series of $u$ with non-negative exponents. Since we do not need this property in this paper, we will not give its proof. Besides all this, it would be interesting
to understand the nature of the function $\bsb{d}_3$. Computer-assisted experiments are possible and
generate large tables of coefficients of the functions $\psi$ and $\bsb{d}_3$, but we will not report them here.

\begin{Theorem}\label{gkstarexpansion}
For all $k\geq 0$, we have the identity:
$$g^\star_k=h^{q^k}\left\{\bsb{d}_2^{(k+1)}\prod_{i=1}^k(t-\theta^{q^i})\bsb{d}_3^\star-\bsb{d}_2\left(\bsb{d}_2^{(k+1)}\left(\frac{1}{u^{q^k}}+\sum_{i=0}^{k-1}\frac{(t-\theta^{q^k})\cdots(t-\theta^{q^{i+1}})}{u^{q^i}}\right)+\frac{(\bsb{d}_3^\star)^{(k+1)}}{t-\theta^{q^{k+1}}}\right)\right\}.$$
\end{Theorem}

\noindent\emph{Proof.} 
First of all, we recall that, for $k\geq 0$,
$$\tau^{k+1}s_\carlitz=L_k^\star s_\carlitz,$$
where $L_k^\star=(t-\theta^{q^k})\cdots(t-\theta)$.
We also recall that $F^{(1)}=F+R$, so that
\begin{eqnarray*}
F^{(k+1)}&=&F+\sum_{i=0}^kR^{(i)}\\
&=&F+\sum_{i=0}^k\frac{1}{L_i^\star s_\carlitz u^{q^i}}\\
&=&s_\carlitz^{-1}\left(F^\star+\sum_{i=0}^k\frac{1}{L_i^\star u^{q^i}}\right).
\end{eqnarray*}

Moreover,
$$\bsb{d}_3=\frac{\bsb{d}_3^\star}{L_0^\star s_\carlitz},$$
yielding
$$\bsb{d}_3^{(k+1)}=\frac{(\bsb{d}_3^\star)^{(k+1)}}{L_{k+1}^\star s_\carlitz}.$$

Therefore, by (\ref{intermd1}) we deduce:
\begin{eqnarray*}
\tau^{k+1}(s_\carlitz\bsb{d}_1)&=&L_k^\star s_\carlitz\bsb{d}_1^{(k+1)}\\
&=&L_k^\star s_\carlitz\left(\bsb{d}_2^{(k+1)}F^{(k+1)}+\bsb{d}_3^{(k+1)}\right)\\
&=&L_k^\star s_\carlitz\left(\bsb{d}_2^{(k+1)}\left(s_\carlitz^{-1}\left(F^\star+\sum_{i=0}^k\frac{1}{L_i^\star u^{q^i}}\right)\right)+\frac{(\bsb{d}_3^\star)^{(k+1)}}{L_{k+1}^\star s_\carlitz}\right)\\
&=&L_k^\star\left(\bsb{d}_2^{(k+1)}\left(F^\star+\sum_{i=0}^k\frac{1}{L_i^\star u^{q^i}}\right)+\frac{(\bsb{d}_3^\star)^{(k+1)}}{L_{k+1}^\star}\right)\\
&=&L_k^\star\bsb{d}_2^{(k+1)}F^\star+\bsb{d}_2^{(k+1)}\sum_{i=0}^k\frac{L_k^\star}{L_i^\star u^{q^i}}+\frac{(\bsb{d}_3^\star)^{(k+1)}}{t-\theta^{q^{k+1}}}.
\end{eqnarray*}
Furthermore,
\begin{eqnarray*}
(\tau^{k+1}s_\carlitz)\bsb{d}_1&=&L_k^\star s_\carlitz \bsb{d}_1\\
&=&L_k^\star s_\carlitz(\bsb{d}_2s_\carlitz^{-1}F^\star+(t-\theta)^{-1}s_\carlitz^{-1}\bsb{d}_3^\star)\\
&=&L_k^\star\bsb{d}_2F^\star+\prod_{i=1}^k(t-\theta^{q^i})\bsb{d}_3^\star.
\end{eqnarray*}
Subtracting, the terms containing $F^\star$ cancel each-others and we obtain the formula applying Theorem \ref{charroots}.\CVD

\noindent\emph{Remark.} Theorem \ref{gkstarexpansion} allows to compute the $u$-expansions
of $(\bsb{d}_3^\star)^{(k)}|_{t=\theta}$ for all $k$. For example, we deduce from the identity of the theorem for $k=0$,
$$\frac{-1}{h}=\frac{\bsb{d}_2\bsb{d}_2^{(1)}}{u}+\frac{1}{t-\theta^q}(\bsb{d}_3^\star)^{(1)}\bsb{d}_2-(\bsb{d}_3^\star)\bsb{d}_2^{(1)},$$
after evaluation at $t=\theta$:
\begin{eqnarray*}
(\bsb{d}_3^\star)^{(1)}|_{t=\theta}&=&\left(\frac{E}{u}-1\right)\frac{\theta-\theta^q}{h}.
\end{eqnarray*}

\noindent\emph{Proof of Corollary \ref{corgstar}.} We assume here that $q\neq2$ (but the case $q=2$ can 
be handled in a similar way, with slightly different results). We compute the truncation $[g_k^\star]_{k+1}$ to the order $q^{k+1}$,
by using the following properties:
$$[\bsb{d}_2^{(k+1)}]_{k+1}=1,\quad [h^{q^k}]_{k+1}=-u^{q^k},\quad [(\bsb{d}_3^\star)^{(k+1)}]_{k+1}=0$$
and we proceed as in the proof of Proposition \ref{lemmed1}. We decompose the sum in the right-hand side of
the formula of Theorem \ref{gkstarexpansion} in four terms:
\begin{eqnarray*}
{[}h^{q^k}\bsb{d}_2\bsb{d}_2^{(k+1)}u^{-q^k}{]}_{k+1}&=&-{[}\bsb{d}_2{]}_{k+1}\\
\left[h^{q^k}\bsb{d}_2\bsb{d}_2^{(k+1)}\sum_{i=0}^{k-1}\frac{(t-\theta^{q^k})\cdots(t-\theta^{q^{i+1}})}{u^{q^i}}\right]_{k+1}&=&
-\left[\bsb{d}_2\sum_{i=0}^{k-1}(t-\theta^{q^k})\cdots(t-\theta^{q^{i+1}})u^{q^k-q^i}\right]_{k+1},\\
\left[h^{q^k}\bsb{d}_2\bsb{d}_2^{(k+1)}\frac{(\bsb{d}_3^\star)^{(k+1)}}{t-\theta^{q^{k+1}}}\right]_{k+1}&=&0,\\
\left[-h^{q^k}\bsb{d}_2^{(k+1)}\prod_{i=1}^{k}(t-\theta^{q^i})\bsb{d}_3^\star\right]_{k+1}&=&
\left[u^{q^k}\bsb{d}_3^\star\right]_{k+1}\prod_{i=1}^{k}(t-\theta^{q^i}).
\end{eqnarray*}
The corollary follows summing up everything and using (\ref{qdifferent2}).\CVD

\subsection{Appendix: transcendence of $F^\star$ and $\bsb{d}_1$ over formal Laurent series}

Although we will not need it in this paper, we prove here, for further references, the transcendence of $F^\star$ and $\bsb{d}_1$
over the field $\CC_\infty(t)((u))$.

By \cite{Ke1}, we can embed an algebraic closure of $\CC_\infty(t)((u))$ in the ring $\CC_\infty(t)^{\text{alg.}}\langle\langle u\rangle\rangle$ of 
{\em generalised formal series} $\sum_{i\in\mathcal{I}}c_iu^i$ (whose support, ordered with $\leq$, is a {\em well ordered subset} of $\QQ$; see definition in loc. cit.). We choose such an embedding.

\begin{Lemme}\label{lemmatranscendence} 
The function $F^\star$ is transcendental over the field 
$\CC_\infty(t)((u))$.
\end{Lemme}

\noindent\emph{Proof.} The function $F^\star$ is identified in an unique way with a generalised formal series.
The functional equation ensures that this series has the following $u$-expansion:
$$F^\star=\sum_{n\geq0}c_nu^{-1/q^{n+1}},$$ for some $c_0,c_1,\ldots$ in $\CC_\infty(t)^{\text{alg.}}$.

Actually, these coefficients can be computed easily, by using the functional equation,
(\ref{valueFtheta}) and the limit $\lim_{t\to\theta}(t-\theta)s_\carlitz(t)=s_\carlitz^{(1)}(\theta)=-\widetilde{\pi}$. Although we will not use them here,
we give their formulas for the sake of completeness:
$c_0=1$ and
$$c_n=(t-\theta)(t-\theta^{1/q})\cdots(t-\theta^{1/q^n}),\quad n>0.$$

Let us suppose by contradiction that $F^\star$ is algebraic over $\CC_\infty(t)((u))$. By \cite[Theorem 8]{Ke1} (read also the discussion on top of page 3465
and \cite{Va}),
there exist $k$ and $d_0,d_1,\ldots,d_k\in \CC_\infty(t)^{\text{alg.}}$, not all zero, such that, for all $n$,
\begin{equation}\label{eqequationLRR}d_0c_n+d_1c_{n+1}^p+\cdots+d_{k}c_{n+k}^{p^k}=0,
\end{equation} where $p$ is the prime dividing $q$.

This means that $F^\star$ is algebraic over $\CC_\infty(t,u)$. Consider now the completion $\mathcal{L}_\infty=\CC_\infty(t)^{\text{alg.}}((u^{-1}))$ of $\CC_\infty(t)^{\text{alg.}}(u)$
for the $u^{-1}$-valuation. Then, the image of $F^\star$ in $\mathcal{L}_\infty$ can be identified 
with a double formal series of $\CC_\infty((t))((u^{-1}))$
which
converges at every $(t,u)$ such that $|u|>1$ and $|t|<q$
to the function $$G:u\mapsto\sum_{i\geq 0}e_\carlitz\left(\frac{\log_{\carlitz}(u^{-1})}{\theta^{i+1}}\right)t^i,$$
where $\log_\carlitz$ is the Carlitz's logarithmic series. 

The latter function extends to the $u$'s such that $|u|>|\widetilde{\pi}|^{-1}$
and the function $z\mapsto F^\star(z)$ factors through $G$. Our assumptions imply that for all $z\in \CC_\infty$,
$F^\star(z)\in\TT_{<q}$ is algebraic over $\CC_\infty(t)$. However, if $z=1$ we find $F^\star(1)=s_\carlitz(t)\in \CC_\infty((t))$,
which is a transcendental function.\CVD

\begin{Corollaire}\label{transcendence_of_d1}
The function $s_\carlitz^{-1}\bsb{d}_1$ is transcendental over $\CC_\infty(t)((u))$.
\end{Corollaire}
\noindent\emph{Proof.} We have, by definition, 
\begin{equation}\label{intermd1}s_\carlitz\bsb{d}_1=\bsb{d}_2F^\star+(t-\theta)^{-1}\bsb{d}_3^\star.\end{equation} We know by Proposition \ref{prarchiv} part 5, that 
$\bsb{d}_2$ belongs to $\FF_q[t,\theta][[u]]$. Moreover, by 
Proposition \ref{lemmed1} part 3, we know that $\bsb{d}_3^\star\in\FF_q[t,\theta][[u]]$.
Finally, by Lemma \ref{lemmatranscendence}, $F^\star$ is transcendental over $\CC_\infty(t)((u))$.\CVD

\end{document}